\documentclass[a4paper,reqno]{amsart}
\usepackage{amsmath}
\usepackage{amsthm}   
\usepackage{amssymb}
\usepackage{epsfig}
\usepackage{psfrag}
\usepackage{graphicx}
\usepackage{caption2}
\usepackage{textcomp}
\usepackage{color}

\newtheorem{theorem}{Theorem}[section]

\newtheorem{corollary}[theorem]{Corollary}

\newtheorem{lemma}[theorem]{Lemma}
\newcommand{\Rz}{\mathbb{R}}
\newcommand{\Nz}{\mathbb{N}}
\newcommand{\argmin}{\text{\,\rm Arg\,min\,}}

\newcommand{\lan}{\left\langle}
\newcommand{\ran}{\right\rangle}
\newcommand{\disp}{\displaystyle}
\newcommand{\haz}{\widehat}

\newcommand{\ove}{\overline}
\newcommand{\epsi}{\varepsilon}
\newcommand{\ddt}{\disp\frac{\text{\rm d}}{\rm dt}}
\newcommand{\dd}{\,{\rm d}}
\newcommand{\Rzn}{\Rz^{3 \times 3}_{\,\text{\rm sym}}}
\newcommand{\Rzd}{\Rz^{3 \times 3}_{\, \text{\rm dev}}}
\newcommand{\GammaD}{\Gamma_{\text{\rm Dir}}}
\newcommand{\GammaT}{\Gamma_{\text{\rm tr}}}
\newcommand{\C}{{\mathbb C}}
\newcommand{\HH}{{\mathbb H}}
\newcommand{\FF}{{\mathcal F}}
\newcommand{\tr}{\text{tr}}
\newcommand{\ep}{p}
\renewcommand{\H}{{\mathcal H}}
\newcommand{\calL}{{\mathcal L}}

\begin{document}

\title[A variational principle for hardening elastoplasticity]{A variational principle \\for hardening elastoplasticity} 

\author{Ulisse Stefanelli}
\address{IMATI - CNR, 
v. Ferrata 1, I-27100 Pavia, Italy.}
\email{ulisse.stefanelli\,@\,imati.cnr.it}\thanks{This research was performed during a visit to the Seminar for Applied Mathematics at ETH Z\"urich and the Institute of Mathematics of the University of Z\"urich under the sponsorship of the STM CNR 2006 Program and the Swiss National Science Foundation. The support and hospitality of both institutions is gratefully acknowledged.}
\urladdr{http://www.imati.cnr.it/ulisse/}


\begin{abstract} We present a variational principle governing the quasistatic evolution of a linearized elastoplastic material. In case of linear hardening, the novel characterization allows to recover and partly extend some known results and proves itself to be especially well-suited for discussing general approximation and convergence issues. In particular, the variational principle is exploited in order to prove in a novel setting the convergence of time and space-time discretizations as well as to provide some possible a posteriori error control.
\end{abstract}

\maketitle

\pagestyle{myheadings}

\section{Introduction}\label{intro}
\setcounter{equation}{0}

The primal initial-boundary value problem of elastoplasticity consists in determining the generalized deformation state of a material subject to external mechanical actions. In particular, starting from some initial state and for a given load and traction, one shall  determine the displacement $u$ of the body from the reference configuration, the inelastic (plastic) part $p$ of its strain, and, possibly, a vector of internal hardening variables $\xi$. In the small deformation regime and within the frame of associative elastoplasticity, the problem is classically formulated in a variational form as that of finding the absolutely continuous trajectory $t \in [0,T] \mapsto y(t) \in Y$ ($Y$ is a Banach space) such that 
\begin{equation}
  \label{eq}
  \partial \psi(y') + A y \ni \ell \quad \text{a.e. in} \ \ (0,T), \ \ y(0)=y_0,
\end{equation}
where $y=(u,p,\xi)$ stands for the vector of unknown fields, $A : Y \to Y^*$ (dual) is linear, continuous, and symmetric, and $\psi:Y \to [0,\infty]$ is the positively $1$-homogeneous and convex dissipation potential ($\partial$ is the classical subdifferential in the sense of Convex Analysis, see below). Moreover, $\ell:[0,T]\to Y^*$ is a given and suitably smooth generalized load (possibly including surface tractions) and $y_0$ represents the initial state. The reader is referred to Section \ref{mechanical_model} for some brief mechanical motivation as well as to the classical monographs by {\sc Duvaut \& Lions} \cite{Duvaut-Lions}, {\sc Han \& Reddy} \cite{Han-Reddy99}, {\sc Lemaitre \& Chaboche} \cite{Lemaitre-Chaboche90}, and {\sc Simo \& Hughes} \cite{Simo-Hughes} for a comprehensive collection of results.

The aim of this paper is that of investigating a global-in-time variational formulation of problem \eqref{eq}. In particular, we shall introduce the functional \linebreak $\FF:W^{1,1}(0,T;Y) \to [0,\infty]$ on trajectories as
$$\FF(y)=\int_0^T \big(\psi(y') + \psi^*(\ell - Ay) - \lan \ell - Ay, y'\ran \big),$$
where $\psi^*$ stands for the conjugate $\psi^*(w)=\sup_{v \in Y}(\lan w, v\ran - \psi(v))$ of $\psi$ and $\lan\cdot, \cdot\ran$ denotes the duality pairing between $Y^*$ and $Y$. The starting point of this analysis relies on the fact that solutions of \eqref{eq} and minimizers of $\FF$ fulfilling the initial condition coincide, namely (see Theorem \ref{characterization})
\begin{gather}
  y \ \ \text{solves \eqref{eq}} \quad \text{iff}\quad \FF(y)=\min \FF = 0 \ \ \text{and} \ \ y(0)=y_0. \label{vc}
\end{gather}

This variational characterization has a clear mechanical interpretation. Indeed, since $\psi$ is positively $1$-homogeneous, its conjugate $\psi^*$ turns out to be the indicator function of the convex set $\partial \psi(0)$. Hence, $\min \FF$ is actually a constrained minimization problem and we have
\begin{eqnarray}
  \FF(y)=0 \quad \text{iff} \quad
 \left\{
 \begin{array}{lll}
\ell - Ay \in \partial \psi(0) \quad \text{a.e. in} \ \ (0,T),\\
\varphi(T,y(T)) + \disp\int_0^T\psi(y') = \varphi(0,y(0))- \disp\int_0^T \lan \ell',y \ran,
  \end{array}
\right.\label{vc2}
\end{eqnarray}
where we have used the notation $(t,y) \mapsto \varphi(t,y)=  \frac12 \lan Ay, y\ran - \lan \ell(t), y\ran$. The first relation above expresses the so-called {\it local stability} \cite{Mielke05} of the trajectory whereas the second is nothing but the energy balance at time $T$. More precisely, $ \varphi(t,y)$ denotes the complementary energy at time $t$ for the state $y$, $\int_0^T \psi(y')$ represents the dissipation of the system on $[0,T]$, and $- \int_0^T \lan \ell',y \ran$ is the work of external actions on $[0,T]$. Hence, minimizing $\FF$ consists in selecting the (only) stable trajectory which conserves the energy. In this regards, the reader is referred to the pioneering papers by {\sc Moreau}~\cite{Moreau68,Moreau70,Moreau71}.

The interest of variational characterization \eqref{vc} of the differential problem \eqref{eq} relies on the possibility of exploiting the general tools from the Calculus of Variations. Some care is however required. Indeed, although $\FF$ is convex and lower semicontinuous with respect to the weak topology of $W^{1,1}(0,T;Y)$, the functional generally fails to be coercive. Moreover, one is not just asked to minimize $\FF$ but also to prove that the minimum is $0$. This considerations suggest that the Direct Method is hardly applicable in order to get solutions to \eqref{eq} via the characterization in~\eqref{vc}.

The first issue of this paper is instead that of showing that the variational principle in \eqref{vc} is particularly well-suited for discussing general approximation issues. Since solutions and minimizers coincide, a natural tool in order to frame an abstract approach to limiting procedures within \eqref{eq} is that of considering the corresponding minimum problems via $ \Gamma$-convergence \cite{DeGiorgi-Franzoni75}. As the value of the functional is directly quantified to be $0$ on the minimizers, what is actually needed here for passing to limits are so-called $\Gamma$-$\liminf$ inequalities only and the latter are generally easily available. We shall specifically focus on the case of linear hardening elastoplasticity and apply the above-mentioned perspective in order to recover in a unified and more transparent frame and partly generalize some convergence results for conformal finite elements  (Thm. \ref{space_approx}), time-discretizations (Thm. \ref{convergence_theta}), and fully-discrete space-time approximations (Thm. \ref{s_t}). In particular, for time-discretization we develop a discrete version of the variational principle \eqref{vc} in the same spirit of the theory of {\it variational integrators} \cite{Marsden-West01} (see Subsection \ref{discrete_principle}). This connection entails also some generalized view at the classical discrete time-schemes (see Subsection \ref{discrete_scheme}).

A second novel point of the present variational approach consists in the possibility of exploiting $\FF$ in order to estimate a posteriori some approximation error. By letting $\FF(y)=0$, we will check that (Cor. \ref{unif_dist})
$$ \max_{[0,T]}\frac12\lan A(y-v),y-v \ran \leq \FF(v)\quad \forall v \in W^{1,1}(0,T;Y), \ \ v(0)=y_0.$$
If $A$ shows some coercivity (which is precisely the case of linearized hardening, see Subsection \ref{A_coerc}), and $v$ is the outcome of some approximation procedure, the estimate above may serve as the basis for some a posteriori estimation procedure, possibly headed to adaptivity (see Subsection \ref{adaptivity}). Let us stress that the latter and \eqref{vc2} entail that the distance of a (stable) trajectory from the solution to \eqref{eq} can be uniformly estimated by means of its energy production along the path.

The variational characterization in \eqref{vc} stems from a reinterpretation in the present elastoplastic frame of the celebrated principle by {\sc Brezis \& Elekand} \cite{Brezis-Ekeland76b,Brezis-Ekeland76} for gradient flows of convex functionals. Since its introduction, the latter principle has continuously attracted attention. In particular, it has been exploited in the direction of proving existence \cite{Rios76,Rios76b,Rios78,Rios79,Auchmuty93,Roubicek00,Ghoussoub-Tzou04} (note that the above-mentioned obstructions to the application of the Direct Methods again appear) and the description of long-time dynamics \cite{Lemaire96}. Moreover, the Brezis-Ekeland approach has been adapted to the case of second order \cite{Mabrouk01,Mabrouk03} and doubly nonlinear equations \cite{be} as well.

One has to mention that, of course, \eqref{vc} is not the only possible global-in-time variational characterization of \eqref{eq}. Besides minimizing the $L^2$ space-time norm of the residual (which might be little interesting since the order of the problem is doubled), one has at least to mention {\sc Visintin} \cite{Visintin01}, where generalized solutions are obtained as minimal elements of a certain partial-order relation on the trajectories,
 and the recent contribution by {\sc Mielke \& Ortiz} \cite{Mielke-Ortiz06} where the functional
\begin{equation}\label{MO}
y \mapsto e^{-T/\epsi} \varphi(T,y(T)) + \int_0^T e^{-t/\epsi}\left(\psi(y') + \disp\frac{1}{\epsi} \varphi(t,y) \right)
\end{equation}
is minimized among trajectories with $y(0)=y_0$. Under extra-smoothness conditions on $\psi$ (not fulfilled in the current frame), the Euler-Lagrange equations of the latter functional are
\begin{gather}
  -\epsi D^2\psi(y')y'' + D \psi (y') + Ay = \ell\nonumber\\
y(0)=y_0, \qquad D\psi(y'(T)) + Ay(T)=\ell(T).\nonumber
\end{gather}
In particular, minimizing the functional in \eqref{MO} consists in performing a suitable elliptic (in time) regularization of the problem. In the specific case of $\psi$ positively $1$-homogeneous, the limit $\epsi \to 0$ can be carried out and the minimizers of the functional in \eqref{MO} are proved to converge to the solution of \eqref{eq}. The latter approach is quite different form that of \eqref{vc}. On the one hand, it is much more general as it naturally applies to the non-smooth case as well (no derivatives of $\psi$ and $\phi$ are involved). On the other hand, it relies on an intermediate and somehow unphysical (not causal) $\epsi$-regularized problem and (besides existence) it seems not directly suited for recovering the full extent of our approximation results for the specific case of problem \eqref{eq}.

\section{Mechanical model}\label{mechanical_model}
\setcounter{equation}{0}

Let us provide the reader with a brief introduction to the mechanical setting under consideration. Our aim is just that of recalling some essential features of the models and well as their variational formulation. In particular, we restrain from reporting here an extensive discussion on associative elastoplasticity as the latter can be easily recovered from the many contributions on the subject. The reader is particularly referred to the mentioned monographs for some comprehensive presentation. 

\subsection{Preliminaries}

We will denote by $ \Rzn $ the space of symmetric $ 3 \times 3$ tensors endowed with the natural scalar product $a : b := \text{tr}(ab)= a_{ij}b_{ij}$ (summation convention). The space $ \Rzn $ is orthogonally decomposed as $\Rzn = \Rzd \oplus \Rz \,1_2$, where $ \Rz \,1_2$ is the subspace spanned by the identity 2-tensor $1_2$ and $ \Rzd $ is the subspace of deviatoric symmetric $ 3 \times 3$ tensors. In particular, for all $ a \in \Rzn$, we have that $ a = a_{\text{\rm dev}} + \tr(a)1_2/3$. 

 We shall assume the reference configuration $ \Omega $ to be a non-empty, bounded, and connected open set in $ \Rz^3 $ with Lipschitz continuous boundary. The space dimension $ 3 $ plays essentially no role throughout the analysis and we would be in the position of reformulating our results in $ \Rz^d $ with no particular intricacy. Our unknown variables are the displacement of the body $u\in  \Rz^3$, the plastic strain $p\in\Rzd$, and a vector of internal variables $\xi \in \Rz^m \ (m \in \Nz)$ which will describe the hardening of the material. We will denote by $ \epsi(u)$ the standard symmetric gradient.

\subsection{Constitutive relation}

Moving within the small-strain regime, we additively decompose the linearized deformation $\epsi(u)$ into the elastic strain $e$ and the inelastic (or plastic) strain $\ep$ as 
$$ \epsi(u) = e + \ep.$$

Let $ \C $ be the elasticity tensor. By regarding the latter as a symmetric positive definite linear map $ \C : \Rzn \rightarrow \Rzn$, we shall assume that the orthogonal subspaces $ \Rzd $ and $ \Rz \, 1_2$ are invariant under $ \C $. This amounts to say that indeed 
$$ \C a= \C_{\text{\rm dev}} a_{\text{\rm dev}} + \kappa\, \text{tr}(a)1_2,$$ 
for a given $ \C_{\text{\rm dev}}: \Rzd \rightarrow \Rzd$ and a constant $ \kappa$, and all $ a \in \Rzn$. The case of isotropic materials is given by $ \C_{\text{\rm dev}} = 2G (1_4 - 1_2 \otimes 1_2/3) $ and $ G $ and $ \kappa $ are respectively the shear and the bulk moduli. The latter decomposition is not exploited in our analysis but it is clearly suggested by the mechanical application.
Moreover, we shall introduce two linear symmetric positive semi-definite hardening moduli $\HH_p:\Rzd\to \Rzd$ and $\HH_\xi:\Rz^m\to \Rz^m$ (to be identified with a fourth order tensor and a matrix, respectively) and define the Helmholtz free energy $W: \Rzn \times \Rzd \times \Rz^m\to [0,\infty)$ of the material as
 $$W(\epsi(u),\ep,\xi):=  \frac12 (\epsi(u) - \ep):\C(\epsi(u) - \ep) + \frac12 \ep : \HH_p \ep + \frac12 \xi^T \cdot \HH_\xi \xi .$$

The generalized stresses $(\sigma,\eta)$ are conjugate to the above-defined generalized strains $(e,\xi)$ via the energy $W$. In particular,
the material is  classically assumed to show elastic response,
\begin{equation}\label{elastic} \sigma = \frac{\partial W}{\partial e}= \C e = \C (\epsi(u) -\ep).
\end{equation}
and the thermodynamic force $\eta$ driving the evolution of the internal variables $\xi$ is defined as
\begin{equation}\label{elastic2}   \eta = - \frac{\partial W}{\partial \xi} = - \HH_\xi \xi.
\end{equation}
Moreover, moving within the frame of associative elastoplasticity, we assume the existence of a function $R: \Rzd \times \Rz^m \to [0,\infty] $ convex, positively $1$-homogeneous, and lower semicontinuous
such that 
\begin{equation}\label{flow}
\partial R(\dot \ep,\dot \xi) \ni \binom{\sigma -\HH_p p}{\eta}.\end{equation}
In particular, $R$ is asked to be the support function of a convex set $C^* \in \Rzd\times \Rz^m$, i.e.  $R(p)= \sup_{q\in C^*} q:p$. We will indicate with $ R^*$ its conjugate, namely the indicator function of $C^*$ given by $R^*(q)= 0$ if $q \in C^*$ and $R^*(q)=\infty$ otherwise. Moreover, we let $C$ be the domain of $R$, namely $C= D(R)=\{(p,\xi) \in \Rzd \times \Rz^m \ : \ R(p,\xi)<\infty\}$.

Finally, the above material relations \eqref{elastic}-\eqref{flow} can be condensed as the following constitutive material law
\begin{equation}\label{const}
 \partial R(\dot \ep, \dot \xi) + \binom{(\C + \HH_p) \ep}{\HH_\xi \xi} \ni  \binom{\C \epsi(u)}{0},
\end{equation}
which in turn can be rephrased in the form of \eqref{eq} by letting
\begin{gather}
 y=(p,\xi),   \   \ Y= \Rzd \times \Rz^m,\ \ \psi=R,\nonumber\\
 A(p,\xi)=\big((\C+\HH_p) \ep,\HH_\xi \xi\big), \ \ \ell = (\C \epsi(u),0).\label{Y_const}
\end{gather}
 Let us close this subsection by explicitly mentioning three classical linear hardening models \cite[Ex. 4.8, p.88]{Han-Reddy99}
\begin{description}
\item[Linear kinetic hardening] choose $\HH_p = h_p 1_4$ where $h_p >0$, and  $\HH_\xi =0$. In this case the internal variable $\xi$ is not evolving and shall be removed from the set of unknowns.
\item[Linear isotropic hardening] choose $\HH_p = 0$, $m=1$, and $\HH_\xi = h_\xi  >0$. Moreover, let $D(R)=\{(p,\xi)\in \Rzd\times \Rz\ : \ |p| \leq \xi\}$.
\item[Linear combined kinetic-isotropic hardening] let $\HH_p = h_p 1_4$, $m=1$, and $\HH_\xi = h_\xi$ where $h_p, h_\xi>0 $. Moreover, let $D(R)=\{(p,\xi)\in \Rzd\times \Rz\ : \ |p| \leq \xi\}$. 
\end{description}
It is  beyond the purpose if this introduction to discuss and justify the above-mentioned material models. The reader shall check the cited references for comments on their relevance within applications and some mechanical motivation.

\subsection{Variational formulation of the quasistatic evolution}

Let us now move to the consideration of the full equilibrium problem. To this aim, we assume that the boundary $ \partial \Omega $ is partitioned in two disjoint open sets $ \GammaT $ and $ \GammaD$ with $ \partial \GammaT = \partial \GammaD$ (in $ \partial \Omega$). We ask $ \GammaD $ to be such that there exists a positive constant $ c_{\text{\rm Korn}} $ depending on $ \GammaD$ and $ \Omega $ such that the Korn inequality
\begin{equation}\label{korn}
c_{\text{\rm Korn}}\| u \|^2_{H^1(\Omega;\Rz^3)} \leq \| u \|^2_{L^2(\GammaD;\Rz^3)} + \| \epsi (u) \|^2_{L^2(\Omega; \Rzn)},
\end{equation}
holds true for all $ u \in H^1(\Omega;\Rz^3)$. It would indeed suffice to impose $ \GammaD $ to have a positive surface measure (see, e.g., \cite[Thm. 3.1, p. 110]{Duvaut-Lions}).

 For the sake of simplicity, we will prescribe homogeneous Dirichlet boundary conditions on $ \GammaD$ (our analysis extends with little notational intricacy to the case of non-homogeneous  Dirichlet boundary conditions as well). On $ \GammaT $ some time-dependent traction will be prescribed instead.

 As for the full quasistatic evolution of the material we shall couple the constitutive relation \eqref{const} with the equilibrium equation
\begin{equation}\label{equi}
\text{div}\,\sigma + f =0 \quad \text{in} \ \ \Omega.
\end{equation}
Here, we assume to be given the body force $ f :[0,T] \to L^2(\Omega;\Rz^3)$ and a surface traction $ g : [0,T] \to L^2(\GammaT;\Rz^3)$. 

Then, one can rephrase the problem into the form of \eqref{eq} by choosing
\begin{eqnarray}
  &&y = (u,\ep,\xi), \label{Y1}\\
&&Y=\Big\{(u,p,\xi) \in H^1(\Omega;\Rz^3)\times L^2(\Omega;\Rzd)\times L^2(\Omega;\Rz^m) \nonumber\\
&&\qquad \qquad \text{such that} \ \  u =0 \ \ \text{on} \ \ \GammaD\Big\},\label{Y2} \\
&&\lan A(u,\ep,\xi), (v,q,z)\ran = \int_\Omega  \Big( (\epsi(u) - \ep): \C(\epsi(v) - q) + \ep : \HH_p q + \xi^T \cdot \HH_\xi z  \Big)\nonumber\\
&&\qquad\qquad \forall (v,q,z) \in Y,\label{Y3} \\
&&\psi(u,\ep,\xi)= \int_\Omega R(\ep,\xi), \label{Y4}
\end{eqnarray}
and defining the total load $ \ell: [0,T] \to Y^*$  as
$$\lan \ell(t), (u,p,\xi) \ran = \int_\Omega f \cdot u + \int_{\GammaT} g \cdot u \, \dd \H^{2} \quad \forall u \in H^1(\Omega; \Rz^3), \ t \in [0,T],$$
where $ \H^2 $ is the 2-dimensional Hausdorff measure.

\subsection{The coercivity of $A$}\label{A_coerc}
Let us close this introductory discussion by explicitly commenting on the coercivity of the bilinear form induced by $A$. We shall recall some sufficient conditions on $\HH_p$, $\HH_\xi$, and $R$ in such a way that there exists a constant $\alpha>0$ such that
\begin{equation}\label{target}
\lan A y , y\ran \geq \alpha |y|^2 \quad \forall y \in D(\psi)
\end{equation}
where $|\cdot|$ is the norm in $Y$. This issue is fairly classical \cite[Sec. 7.3, p. 167]{Han-Reddy99} and we discuss it here for the sake of completeness only. 

Of course \eqref{target} holds (and even for all $y \in Y$) whenever $\HH_p$ and $\HH_\xi$ are positive definite (this is the case of the above-mentioned linear combined kinematic-isotropic hardening).

As we have already observed, in case $\HH_\xi=0$, the problem naturally reduces to the pair $(u,p)$ only. Up to this reduction, \eqref{target} holds (again for all $y \in Y$) when $\HH_p$ is positive definite. This is exactly the case of linear kinematic hardening. 

On the other hand, in case $\HH_p=0$, the plastic strain will still evolve and one has \eqref{target} if $D(R)$ is bounded in the $p$-direction for all $\xi$, namely if \cite[(7.51)]{Han-Reddy99}
\begin{equation}\label{D_bound} D(R)\subset \{(p,\xi) \in \Rzd \times \Rz^m\ : \ \  \beta |p|^2 \leq \xi^T \cdot \HH_\xi \xi \ \ \text{for some constant} \ \ \beta >0\},
\end{equation}
which is clearly the case for linear isotropic hardening. 

Some generalization of the latter condition could in principle be considered for the case when $\HH_p$ and $\HH_\xi$ are only semi-definite. In particular, \eqref{target} holds if one assumes \eqref{D_bound} and
$$ \xi \not =0  \ \ \text{and} \ \ \xi^T \cdot \HH_\xi \xi =0  \ \ \Rightarrow  \ \ R(p,\xi)=\infty \ \ \forall p \in \Rzd.$$

Let us mention that the most critical case in the class of \eqref{const} is $\HH_p=0$, $\HH_\xi=0$ where actually no hardening takes place. This is the situation {\it perfect plasticity} for which the Sobolev space framework above is not appropriate and one would consider the space $BD(\Omega)$ of functions of bounded deformations instead \cite{DalMaso05}. We shall make clear that, even if our variational characterization covers the case of perfect plasticity, the subsequent approximation results apply to the linear hardening situation only.

\section{Characterization}\label{section_characterization}
\setcounter{equation}{0}

\subsection{General assumptions}\label{assumptions}

Let us start by recalling notation and enlisting the basic assumptions for the following. First of all, we will ask that
\begin{equation}
  \label{Y}
  Y \quad \text{is a separable and reflexive Banach space.}
\end{equation}
We will use the symbols $|\cdot|$ for the norm of $Y$ and $\lan \cdot, \cdot \ran$ for the duality pairing between $Y^*$ (dual) and $Y$. The norm in $Y^*$ will be denoted by $|\cdot|_*$ instead.

 We introduce the functional
\begin{gather}
  \psi: Y \to [0,\infty] \quad \text{proper, convex, lower semicontinuous,}\nonumber\\
\text{and positively $1$-homogeneous.} \label{psi}
\end{gather}
Equivalently, $ \psi$ is required to be the support function of a convex and closed set $ C^* \subset Y^* $ containing $ 0$, namely 
$$\psi(y)= \sup\{\lan y^*,y\ran \ : \ y^* \in C^*\}.$$
We shall define $ C=D(\psi)$. Hence, the conjugate $\psi^*: Y^* \to [0,\infty]$, which is classically defined as $\psi^*(y^*)=\sup_{y\in Y}(\lan y^*, y \ran - \psi(y))$, is the indicator function of the convex set $C^*$, namely $\psi^*(y^*)=0$ if $y^* \in C^*$ and $\psi^*(y^*)=\infty$ otherwise. Let us remark that $\psi$ fulfills the triangle inequality $\psi(a) \leq \psi(b) + \psi(c)$ whenever $a= b+c$. 

We shall use the symbol $\partial$ in order to denote the usual subdifferential in the sense of Convex Analysis, namely
$$y^* \in \partial \psi(y) \quad \text{iff}\quad y \in D(\psi) \ \ \text{and} \ \ \lan y^*, w- y\ran \leq \psi(w) - \psi(y) \quad \forall w \in Y.$$
Similarly, we define 
\begin{gather}
y \in \partial\psi^*(y^*)\quad \text{iff}\quad y^* \in D(\psi^*) \ \ \text{and} \ \ \lan w^*- y^*, y\ran \leq \psi^*(w^*) - \psi^*(y^*) \quad \forall w^* \in Y^*\nonumber\\
\quad \text{iff}\quad y^* \in C^* \ \ \text{and} \ \ \lan w^*- y^*, y\ran \leq 0 \quad \forall w^* \in C^*.\nonumber
\end{gather}
Finally, we recall Fenchel's inequality 
$$\psi(y) + \psi^*(y^*) \geq \lan y^*, y \ran\quad \forall y\in Y, \ y^* \in Y^*,$$
and remark that equality holds iff $y^* \in \partial \psi(y)$ (or, equivalently, $y \in \partial^*\psi^*(y^*)$).

As for the operator $A$ we require
\begin{equation}
  \label{A}
  A : Y \to Y^*\quad \text{linear, continuous, and symmetric,}
\end{equation}
and define the function
$$y \to  \phi(y)= \frac12 \lan Ay,y \ran, $$
so that $\, A = D\phi$.
Moreover, we will ask $\phi$ to be coercive on $C=D(\psi)$, namely we assume that there exists a positive constant $\alpha$ such that 
\begin{equation}\label{coerc1}
\phi(y) \geq \frac{\alpha}{2}|y|^2
\quad \forall y \in C.
\end{equation}
As we have already commented in Subsection \ref{A_coerc}, the latter coercivity is fulfilled in the situation of elastoplastic evolution with linear kinematic, isotropic, or combined kinematic-isotropic hardening and will turn out to be sufficient for both the forthcoming characterization results. 

On the other hand, the following uniqueness-type results will be checked under some stronger coercivity frame and we will ask for
\begin{equation}\label{coerc2}
\phi(y) \geq \frac{\alpha}{2}|y|^2
\quad \forall y \in C-C.
\end{equation}
Clearly, condition \eqref{coerc2} is fulfilled when $\phi$ happens to be coercive on the whole space $Y$. The latter applies in particular to the case of linear kinematic and combined kinematic-isotropic hardening elastoplasticity. In this case, $\phi$ defines an equivalent (squared) norm in $Y$.

We shall make use of the following notation 
$$\chi(y)= \phi(y)+ |y|^2 \quad \forall y \in Y$$
Indeed the latter choice is just motivated by simplicity and could be replaced as well by any other $\chi: Y \to [0,\infty)$ such that $\chi(y)=0$ iff $y=0$ and $y \mapsto \chi(y)-\phi(y)$ is lower semicontinuous.

Finally, we shall fix data such that 
\begin{equation}
  \label{data}
  \ell \in L^\infty(0,T;Y^*), \ y_0 \in C.
\end{equation}
The restriction on the choice of the initial datum in $C$ is motivated by the coercivity assumption on $\phi$ in \eqref{coerc1}. On the other hand, we shall explicitly mention that the usual choice for $y_0$ in elastoplasticity is $y_0=0$.

In the forthcoming of the paper the above assumptions \eqref{Y}-\eqref{coerc1} and \eqref{data} will be tacitly assumed (unless explicitly stated). It should be however clear that the above choice is motivated by the sake of simplicity. Indeed, most of the following results still hold under suitably weaker assumptions, as we shall comment. 

\subsection{The functional}\label{functional}

Let the {\it Lagrangian} $ L : (0,T) \times Y \times Y \to [0,\infty]$ be defined as
\begin{gather}
  L(t,y,p)= \psi(p) + \psi^*(\ell(t) - A y) - \lan \ell(t) - Ay ,p\ran \nonumber\\
 \text{for a.e.} \ \ t \in (0,T), \  \forall y,\, p \in Y, \label{L}
\end{gather}
and the functional $\, F:W^{1,1}(0,T;Y) \to [0,\infty]$ as
\begin{equation}
  \label{F}
  F(y)= \int_0^T L(t,y(t),y'(t))\, \dd t + \chi(y(0) - y_0).
\end{equation}

Now, by simply using the chain rule, we obtain that
$$F(y)=\int_0^T \Big( \psi(y') + \psi^*(\ell - A y) - \lan \ell, y'\ran \Big) + \phi(y(T)) - \phi(y(0)) + \chi(y(0) -y_0).$$
A first remark is that, by exploiting the particular form of $\chi$,
\begin{eqnarray}
F(y)&=&\int_0^T \Big( \psi(y') + \psi^*(\ell - A y) - \lan \ell, y'\ran \Big) \nonumber\\
&+& \phi(y(T)) + \phi(y_0) - \lan A y(0) ,y_0\ran + |y(0) - y_0|^2.
\end{eqnarray}
In particular, $F$ is clearly convex.

\subsection{The characterization}\label{subsec_characterization}

\begin{theorem}[Variational principle]\label{characterization}
 $y$ solves \eqref{eq} iff $F(y)=0=\min F$.
\end{theorem}

\begin{proof}
  Owing to Fenchel's inequality we have that
$$L(t,y,p) = 0 \quad \text{iff} \ \ \ell(t) - A y \in \partial \psi(p)$$
and, clearly, $\chi(y(0)-y_0)=0$ iff $y(0)=y_0$.
Hence, all solutions $y$ of \eqref{eq} are such that $F(y)=0$ and vice-versa.
\end{proof}

Let us remark that the latter variational characterization result holds in much grater generality. The proof made no use of the separability and reflexivity of $Y$ nor of the linearity of $A$ (besides its being single-valued and such that $t \mapsto Ay(t)$ is measurable). Moreover, the positive $1$-homogeneity of $\psi$ is unessential \cite{owr}. In particular, the variational approach of Theorem \ref{characterization} can be directly extended to a variety of different dissipative systems possibly including viscous evolution as well. We shall address this perspective in a forthcoming contribution.

We have already observed that $F$ is convex. Moreover, $F$ is lower semicontinuous with respect to the weak topology of $W^{1,1}(0,T;Y)$ since all weakly convergent sequences in $W^{1,1}(0,T;Y)$ are pointwise weakly convergent as well. Hence, one could be tempted to use the Direct Method in order to get the existence of minimizers, i.e. solutions to equation \eqref{eq}. As we commented in the Introduction, this seems to be no trivial task.

First of all, the functional $F$ need not be coercive with respect to the weak topology of $W^{1,1}(0,T;Y)$. Indeed, the functional $\psi$ may degenerate and hence not control the norm of its argument. Moreover, even in the case when $\psi$ is non-degenerate, the homogeneity assumption just entails that the sublevels of $F$ are bounded in $W^{1,1}(0,T;Y)$ and no weak compactness follows. 

Secondly, even assuming coercivity in the weak topology of $W^{1,1}(0,T;Y)$, one would still need to prove that the minimum $ 0$ is attained. This very obstruction to the use of the Direct Method occurs for the Brezis \& Ekeland principle for gradient flows \cite[Rem. 1]{Brezis-Ekeland76} and for its doubly nonlinear generalization in \cite{be} as well.

\subsection{The variational principle for hardening elastoplasticity}\label{variational_plast}
By referring to the notations of Section \ref{mechanical_model}, let us now present the actual form of the functional $F$ for the case of the constitutive relation for linearized elastoplastic materials with linear hardening (see \eqref{Y_const}). In this case the functional reads
\begin{eqnarray}
 F(p,\xi)&=&
\disp\int_0^T \Big( R(\dot p, \dot \xi) + R^*\big(\C (\epsi(u) -p)- \HH_p p,-\HH_\xi \xi\big)\Big)\nonumber\\
&-& \int_0^T \Big(\big(\C (\epsi(u)-p) - \HH_p p \big): \dot p - \xi^T \cdot \HH_\xi \dot \xi \Big)\nonumber\\
&+& \frac12  (\ep(0)-{\ep}_0): (\C+\HH_p)(\ep(0)-{\ep}_0) + \frac12 (\xi(0) - \xi_0)^T\cdot \HH_\xi (\xi(0) - \xi_0) \nonumber\\
&+& |({\ep}(0),\xi(0)) - ({\ep}_0,\xi_0)|^2, \nonumber
\end{eqnarray}
for some given initial datum $({\ep}_0,\xi_0)\in \Rzd \times \Rz^m$ and $\epsi(u) \in L^\infty(0,T;\Rzn)$. 

In the situation of the quasistatic evolution, for some given initial datum \linebreak$(u_0,{\ep}_0, \xi_0)\in Y$, a load $f\in L^\infty(0,T;(L^2(\Omega;\Rz^3)))$, and a traction \linebreak$g\in L^\infty(0,T;(L^2(\GammaT;\Rz^3)))$, the functional reads (see \eqref{Y1}-\eqref{Y4})
\begin{eqnarray}
 F(u,\ep,\xi)&=&
\disp\int_0^T\int_\Omega \Big( R(\dot p, \dot \xi) + R^*\big(\C (\epsi(u) -p)- \HH_p p,-\HH_\xi \xi\big) \Big)\nonumber\\
&-&\disp\int_0^T\int_\Omega f \cdot \dot u- \int_0^T \int_{\GammaT} g \cdot \dot u \dd\H^2\nonumber\\
 &+&\disp\int_0^T \int_\Omega  (\epsi(u) - \ep): \C(\epsi(\dot u) - \dot p) + \ep : \HH_p \dot p + \xi^T \cdot \HH_\xi \dot \xi  \Big)\nonumber  \\
&+&\disp\frac12 \int_\Omega  (\epsi(u(0) - u_0) - (\ep(0) - {\ep}_0)):\C(\epsi(u(0) - u_0) - (\ep(0) - {\ep}_0))\nonumber\\
&+&\disp\frac12\int_\Omega \Big( (\ep(0)-{\ep}_0): \HH_p(\ep(0)-{\ep}_0) +  (\xi(0) - \xi_0)^T\cdot \HH_\xi (\xi(0) - \xi_0)\Big)\nonumber\\
&+&\disp\frac12 \int_\Omega|(u(0),{\ep}(0),\xi(0)) - (u_0,{\ep}_0,\xi_0)|^2,\nonumber
\end{eqnarray}
for all points $(u,p,\xi) \in Y$ such that
\begin{gather}
\disp\int_\Omega (\epsi(u)-p):\C \epsi(v)= \disp\int_\Omega f \cdot v-  \int_{\GammaT} g \cdot v \dd \H^2 \nonumber\\
 \forall v \in H^1(\Omega;\Rz^3)\ \ \text{with} \ \ v=0 \ \ \text{on} \  \ \GammaD, \ \ \text{a.e. in} \ \ (0,T),\nonumber
\end{gather}
and $F(u,p,\xi)=\infty$ otherwise.

\section{Properties of the minimizers}\label{properties}
\setcounter{equation}{0}
For the sake of illustrating the variational principle of Theorem \ref{characterization}, we shall collect here some properties of the trajectories belonging to the domain of the functional $F$ and, in particular, of the minimizers.

\subsection{Trajectories are in $C$}\label{trajectories}

\begin{lemma}\label{C}
Let $F(y)<\infty$. Then $y(t) \in C$ for all $t \in [0,T]$.
\end{lemma}
\begin{proof}
Since $F(y)<\infty$ we have that $y' \in C$ almost everywhere in $(0,T)$. Hence, for all $t\in [0,T]$, we have that $\int_0^t y' \in C$ by Jensen's inequality.
On the other hand $y_0\in C$ and $y(t)=y_0 + \int_0^t y' $. The assertion follows by recalling  that $C$ is a cone.
\end{proof}
\subsection{Stability at regular points of $\ell$}\label{sec_stability}
Assume $\ell : [0,T] \to Y^*$ is given and let the set of {\it stable} states $\, S(t)\subset Y \,$ for $\, t \in [0,T]\,$ be defined as
$$S(t)= \{y \in Y \ : \ \phi(y) - \lan \ell(t), y\ran  \leq \phi(w) - \lan \ell(t), w\ran + \psi(w-y) \ \ \forall w \in Y\}.$$
One can easily prove that 
$$y \in S(t) \quad \text{iff} \quad \ell(t) - Ay \in \partial \phi(0)= C^*.$$

\begin{lemma}[Stability of the minimizers]\label{Stable}
Let $\ell$ be either left- or right-weakly continuous at some point $t \in [0,T]$ and $F(y)< \infty$. Then $y(t)\in S(t)$.
\end{lemma}
\begin{proof}
 Since $F(y)< \infty $, we have that $\, \ell(t) - Ay(t) \in=C^* \,$ for all $ t \in (0,T)\setminus N$, where $|N|=0$. Choose a sequence $\, t_k \in (0,T)\setminus N$ such that $t_k \to t$ (from the left or from the right) and $\ell(t_k) \to \ell(t)$ weakly in $Y^*$. Hence, $\ell(t_k) - A y(t_k)\to \ell(t)- A y(t)$ weakly in $Y^*$ and $\ell(t)- A y(t)\in C^*$.
\end{proof}

In particular, if $\ell$ happens to be right-continuous at $0$, the functional $F$ will not attain the minimum value $0$ unless the initial datum $y_0$ is stable, namely $y_0 \in S(0)$.

\subsection{Equivalent formulations}\label{equivalent}

Letting now $\ell \in W^{1,1}(0,T;Y^*)$, problem \eqref{eq} admits some alternative equivalent formulations \cite[Sec. 2.1]{Mielke05}. We explicitly mention that $y \in W^{1,1}(0,T;Y)$ is said to be an {\it energetic solution} if it solves the {\it energetic formulation} \cite{Mielke-Theil04} of \eqref{eq}, namely
\begin{eqnarray}
&&y(t) \in S(t) \quad \forall t \in [0,T], \label{en_stab}\\
&&\phi(y(t)) - \lan \ell(t) ,y(t) \ran +\int_0^t \psi(y') = \phi(y(0)) - \lan \ell(0) ,y(0) \ran - \int_0^t \lan \ell',y \ran  \nonumber\\
&&\qquad\quad \forall t \in [0,T], \label{en_en}\\
&&y(0)=y_0. \label{en_initial}
\end{eqnarray}

{\sc Mielke \& Theil} \cite{Mielke-Theil04} proved that the latter is equivalent to \eqref{eq} and hence, owing to the characterization of Theorem \ref{characterization}, to $F(y)=0=\min F$. For the aim of pointing out some features of our variational approach, we shall present here a direct proof of this fact.

\begin{lemma}[Equivalence with the energetic formulation] Let $\ell \in W^{1,1}(0,T;Y^*)$. Then,
  $F(y)=0=\min F$ iff $y$ fulfills \eqref{en_stab}-\eqref{en_initial}. 
\end{lemma}

\begin{proof}
  Owing to Lemma \ref{Stable}, we readily have that the stability condition \eqref{en_stab} holds iff $ \psi^*(\ell - Ay)=0$ almost everywhere.

Let $y$ be such that $F(y)=0$. Then \eqref{en_stab} and \eqref{en_initial} hold and $L(t,y(t),y'(t))=0$ for a.e. $t \in (0,T)$. In particular, for all $t\in [0,T]$,
\begin{gather}
0= \int_0^t L(s,y(s),y'(s))\, \dd s= \int_0^t\Big(\psi(y') - \lan \ell,y' \ran\Big) + \phi(y)\Big|_0^t \nonumber\\
= (\phi(y) - \lan \ell, y \ran)\Big|_0^t + \int_0^t\psi(y') + \int_0^t \lan \ell',y \ran, \nonumber
\end{gather}
so that the energy equality \eqref{en_en} holds for all $t \in [0,T]$.

On the contrary, let $y\in W^{1,1}(0,T;Y)$ fulfill \eqref{en_stab}-\eqref{en_initial}. Then $\chi(y(0) - y_0)=0$ and $\psi^*(\ell - Ay) =0$ almost everywhere (see above). Hence $F(y)=0$ follows from the energy equality  \eqref{en_en} at time $T$ and an integration by parts.
\end{proof}

Let us mention that the latter lemma proves in particular that the energy equality \eqref{en_en} could be equivalently enforced {\it at the final time} $T$ only. Moreover, it proves that, as already commented in the Introduction, all {\em stable trajectories} $t \mapsto y(t)$ (i.e. trajectories such that $y(t) \in S(t)$ for all $t \in [0,T]$) are such that the following energy inequality holds
$$\phi(y(t)) - \lan \ell(t) ,y(t) \ran +\int_0^t \psi(y') \geq \phi(y(0)) - \lan \ell(0) ,y(0) \ran - \int_0^t \lan \ell',y \ran  \quad \forall t \in [0,T].$$
Hence, we have provided a proof to \cite[Prop. 5.7]{Mielke05} (in a somehow simpler situation though).

Before closing this subsection, let us explicitly remark that the above inferred equivalence between formulations has been obtained for the absolutely continuous case $\ell \in W^{1,1}(0,T;Y^*)$ only, whereas the characterization of Theorem \ref{characterization} holds more generally for bounded $\ell$.

\subsection{The functional controls the uniform distance: uniqueness}\label{sec_uniqueness}

So far, we have simply reformulated known results in a variational fashion. Here, we present some novel results instead.

\begin{lemma}[Uniform distance control via $F$]\label{uniform_distance} We have that
\begin{gather}
\eta(1-\eta)\max_{t \in [0,T]}\phi\big(u(t)- v(t)\big) \leq\eta F(u) + (1-\eta)F(v)\nonumber\\
 \forall u, \, v \in   W^{1,1}(0,T;Y), \ \eta \in [0,1].\label{star}
\end{gather}
\end{lemma}

\begin{proof}
The statement follows from the quadratic character of $\phi$. Fix $t \in [0,T]$ and define $F^t: W^{1,1}(0,t;Y) \to [0,\infty]$ as
\begin{eqnarray}
F^t(y)&=& \int_0^t L(s,y(s),y'(s))\, \dd s + \chi(y(0) - y_0)\nonumber\\
&=& \int_0^t \Big( \psi(y') + \psi^*(\ell - A y) - \lan \ell, y'\ran \Big)\nonumber\\
 &+& \phi(y(t)) + \phi(y_0) - \lan A y(0) ,y_0\ran + |y(0) - y_0|^2.\nonumber
\end{eqnarray}
Then, clearly $y \mapsto G^t(y)=F^t(y) - \phi(y(t))$ is convex. Hence, letting $w=\eta u + (1-\eta)v$ we have that 
\begin{gather}
  0 \leq F^t(w) \nonumber\\
\leq \eta \big(G^t(u) + \phi(u(t))\big) + (1-\eta)\big(G^t(v) + \phi(v(t))\big) - \eta(1-\eta)\phi\big(u(t)- v(t)\big), \nonumber
\end{gather}
whence the assertion follows.
\end{proof}

The latter lemma exploits the quadratic character of $\phi$ only. In particular, no coercivity for $\phi$ is assumed. It should however be clear that its application (as the title of the lemma indeed suggests) will always be referred to the situation where the stronger \eqref{coerc2} is required, namely when the left-hand side of \eqref{star} controls 
$$\eta(1-\eta)\frac{\alpha}{2}\max_{[0,T]}|u-v|^2.$$ 

We now present two immediate corollaries of Lemma \ref{uniform_distance}.

\begin{corollary}[Uniform distance from the minimizer]\label{unif_dist}
   Let $F(y)=0$. Then
\begin{gather}
\max_{t \in [0,T]}\phi\big(y(t)- v(t)\big) \leq F(v).\nonumber
\end{gather}
\end{corollary}

This corollary encodes an interesting novel feature of our variational approach, for it provides possible a posteriori error estimator to be used within approximation procedures. It is interesting to remark that the uniform distance of any stable trajectory from the minimizer is controlled by means of its energy production along the path only. Again, although Corollary \ref{unif_dist} holds under no coercivity assumptions of $\phi$, let us mention that its application will be restricted to the frame of \eqref{coerc2}. Finally, we have uniqueness of the minimizers of $F$ attaining the value $0$.

\begin{corollary}[Uniqueness]\label{uniqueness} Assume \eqref{coerc2}. Then, there exists at most one trajectory $y$ such that $F(y)=0$.  
\end{corollary}

\subsection{Lipschitz bound}\label{lipschitz_bound}

Form this point on and throughout the remainder of the paper we shall tacitly assume
\begin{equation}
\ell \in W^{1,\infty}(0,T;Y^*), \quad y_0 \in S(0).\label{data2}
\end{equation}

As already commented after Lemma \ref{Stable}, the above restriction on the initial datum is mandatory whenever $\ell$ admits a weak-right-limit in $0$.

As for $\ell$, the extra Lipschitz continuity assumption is motivated by the rate-independence of the problem (every absolutely continuous datum can be time-rescaled to a Lipschitz continuous datum) and the following well-known result. 

\begin{lemma}[Lipschitz bound]\label{lipschitz}
   Assume \eqref{coerc2} and let  $\ell \in W^{1,\infty}(0,T;Y^*)$, and $F(y)=0$. Then 
   \begin{equation}
     \label{a_priori_bound}
     \|y'\|_{L^\infty(0,T;Y)} \leq \frac{1}{\alpha} \|\ell'\|_{L^\infty(0,T;Y^*)} \quad \text{a.e. in} \ \ (0,T).
    \end{equation}
\end{lemma}

The proof of the lemma is exactly the classical one \cite[Thm. 7.5]{Mielke-Theil04}, but formulated by means of our variational arguments. We provide it for the sake of completeness.

\begin{proof}
  Let $0 \leq s < t \leq T$ be fixed. Since $L(y,y')=0$ almost everywhere we have that 
$$\int_s^t \Big(\psi(y') + \lan \ell', y \ran \Big) + (\phi(y) - \lan \ell, y \ran)\Big|_s^t =0.$$
On the other hand, owing to the strong monotonicity of $A$ and the fact that $y(s)\in S(s)$ (see Lemma \ref{Stable}) one obtains
\begin{gather}
  \phi(y(t) - y(s)) \leq \phi(y(t)) - \lan \ell(s),y(t)\ran + \psi(y(t)- y(s)) -  \phi(y(s)) - \lan \ell(s),y(s)\ran\nonumber\\
 =  (\phi(y) - \lan \ell, u \ran)\Big|_s^t + \int_s^t \lan\ell'(r),y(t) \ran\dd r +  \psi(y(t)- y(s)).\nonumber
\end{gather}
By taking the sum of these two relations and recalling that, by Jensen,
$$\psi(y(t)- y(s))= \psi\left(\int_s^t y' \right)\leq \int_s^t \psi(y')$$
we get that
$$\phi(y(t) - y(s)) \leq \int_s^t \lan\ell'(r),y(t) - y(r) \ran\dd r.$$
Finally, an application of some extended Gronwall lemma (see \cite[Thm. 3.4]{Mielke05}) entails that 
$$\frac{\alpha}{2} |y(t) - y(s)|^2 \leq \frac12 \|\ell'\|_{L^\infty(0,T;Y^*)}|y(t) - y(s)|\, (t-s)$$and the assertion follows.

\end{proof}

\section{Space and data approximation}\label{approximation}
\setcounter{equation}{0}

We now apply the characterization results of Theorem \ref{characterization} to the approximation of solutions of \eqref{eq}. As already commented in the Introduction, we shall proceed via $ \Gamma$-convergence \cite{DeGiorgi-Franzoni75}. The reader is referred to the monographs by {\sc Attouch} \cite{Attouch} and {\sc Dal Maso} \cite{DalMaso93} for some comprehensive discussion on this topic. 
Indeed, since Theorem \ref{characterization} directly quantifies the value of the minimum to be $ 0$, what is actually needed for passing to limits are $ \Gamma-\liminf\,$ inequalities only. We shall illustrate this fact by discussing the simple case of data approximation first.

\begin{lemma}[Data approximation]\label{data_approx} 
Assume \eqref{coerc2}, let $\ell_h \to \ell$ strongly in \linebreak$L^1(0,T;Y^*)$ being uniformly Lipschitz continuous, and $y_{0,h} \to y_0$. Moreover, let $F_h:W^{1,1}(0,T;Y) \to [0,\infty]$ be defined as
\begin{gather}
F_h(y)=\int_0^T \Big( \psi(y') + \psi^*(\ell_h - A y) - \lan \ell_h -Ay, y'\ran \Big) + \chi(y(0) -y_{0,h}),\nonumber
\end{gather}
and let $F_h(y_h) =0$. Then $y_h \to y$ weakly star in $W^{1,\infty}(0,T;Y)$ and $F(y)=0$.
\end{lemma}

\begin{proof}
  Owing to Lemma \ref{lipschitz}, we find a (not relabeled) subsequence $y_h$ such that $y_h \to y$ weakly star in $W^{1,\infty}(0,T;Y)$. Hence, we have by lower semicontinuity  \begin{eqnarray}
    0 \leq F(y) &\leq& \liminf_{h \to 0}\Bigg(\int_0^T \Big( \psi(y'_h) + \psi^*(\ell_h - A y_h) - \lan \ell_h, y'_h\ran \Big)\nonumber\\
&+& \phi(y_h(T))+ \phi(y_{0,n}) - \lan A y_h(0),y_{0,h}\ran + |y_h(0) -y_{0,h}|^2\Bigg)\nonumber\\
&=& \liminf_{h \to 0} F_h(y_h) =0.\nonumber
  \end{eqnarray}
Hence, $F(y)=0$, $y$ is unique, and the assertion follows from the fact that the whole sequence converges.
\end{proof}

\subsection{Preliminaries on functional convergence}\label{preliminaries_functional}

In order to move to more general approximation situations, we are forced to discuss a suitable functional convergence notion. We limit ourselves in introducing the relevant definitions, referring to the mentioned monographs for all the necessary details.

Recall that $Y$ is a real reflexive Banach space. Letting $f_n ,\, f:Y \to  (-\infty,\infty] $ be convex, proper, and lower semicontinuous, we say that $f_n \to f $ in the {\em Mosco sense in $Y$} \cite{Attouch,Mosco69} iff,
for all $y \in Y$,
\begin{gather}
  f(y) \leq \liminf\limits_{n \to \infty}f_n(y_n) \quad \forall  y_n \to y \ \ \text{weakly in} \ \ Y,\nonumber\\
\exists y_n \to y  \ \ \text{strongly in} \ \ Y \ \ \text{such that}  \ \ f(y)= \limsup\limits_{n \to \infty}f_n(y_n).\nonumber
\end{gather}
In particular, $ f_n \to f $ in the {Mosco sense} iff $ f_n \to f$ in the sense of $\Gamma$-convergence with respect to both the weak and the strong topology  in $Y $.

We will consider the situation of approximating functionals $\psi_h$. By \cite[Thm. 3.18, p. 295]{Attouch}, we have that  $\psi_h \to \psi $ in the Mosco sense in $Y$ iff $\psi_h^* \to \psi^* $ in the Mosco sense in $Y^*$. By assuming the functionals $\psi_h$ to be positively $1$-homogeneous, it turns out that the Mosco convergence $\psi_h \to \psi $ in $Y$ is equivalent to  the Mosco convergence of sets ${C^*_h} \to {C^*} $ in $Y^*$ which reads, 
\begin{eqnarray}
   &&C^*_n \ni y^*_n \to y^* \quad \text{weakly in $Y^*$} \ \ \Rightarrow \ \ y^* \in C^*,\nonumber\\
&&\forall y^* \in C^*, \ \ \exists  y^*_n\in C^*_n  \ : \ y^*_n \to y^* \quad \text{strongly in $Y^*$}. \nonumber
\end{eqnarray}

Finally, we repeatedly use a lemma from \cite{be} which we report it here, for the sake of completeness.

\begin{lemma}[Cor. 4.4, \cite{be}]\label{cor44} 
Let $p \in [1,\infty]$ and $f_h,\, f: Y \to (-\infty,\infty]$ be convex, proper, and lower semicontinuous such that 
$$f(y) \leq \inf\left\{\liminf\limits_{h \to 0}f_h(y_h) \ : \ y_h \to y \ \ \text{weakly in} \ \ Y \right\}\quad \forall y \in Y.$$
Moreover, let $y_h \to y$ weakly in $W^{1,p}(0,T;Y)$ (weakly star if $p=\infty$), Then, we have that
$$\int_0^Tf(y(t))\,\dd t \leq \liminf\limits_{h \to 0}\int_0^T f_h(y_h(t))\,\dd t.$$
\end{lemma}

\subsection{Space approximations}\label{approx_space}

We now move to the analysis of some  space approximation situation, indeed specifically tailored for the case of conformal finite elements. Let us enlist here our assumptions for the sake of later referencing.

We assume to be given
\begin{eqnarray}
&& Y_h \subset Y \ \ \text{closed subspaces such that} \ \ \bigcup_{h>0} Y_h \ \ \text{is dense in} \ \ Y,\label{Y_h}\\
  &&\phi_h(y) = \phi(y) \ \ \text{if} \ \ y \in Y_h \ \ \text{and} \ \ \phi_h(y)=\infty \ \ \text{otherwise}.\label{phi_h}\\
 &&\psi_h: Y \to (-\infty,\infty] \quad \text{convex, proper, and lower semicontinuous}, \label{psi_h}\\
 &&\psi_h  \quad \text{positively $1$-homogeneous}, \label{psi_h2}\\
 &&\psi_h \to \psi \quad \text{in the Mosco sense in $Y$}, \label{psi_h3}\\
&&\phi(y) \geq \frac{\alpha}{2}|y|^2
 \quad \forall y \in C_h - C_h \quad \text{where} \ \ C_h=D(\psi_h),\label{phi_psi}\\
 &&\ell_h \to \ell \quad \text{pointwise strongly in $Y^*$},\label{ell_h}\\
 &&\ell_h \quad \text{uniformly Lipschitz continuous},\label{ell_h2}\\
&&y_{0,h}\in Y_h, \ \ y_{0,h} \to y_0.\label{y_0h}
\end{eqnarray}

We shall mention that within the frame of conformal finite elements methods the subspaces $Y_h$ are obviously taken to be finite-dimensional and that the approximating functionals $\phi_h$ and $\psi_h$ are usually the restrictions of the functionals $\phi$ and $\psi$ on the subspace $Y_h$. This is exactly our choice here for $\phi_h$. In particular, one shall observe that $\phi_h \to \phi$ in the Mosco sense in $Y$, $D(\partial \phi_h) = Y_h$, and that
\begin{equation}\label{identity_A} 
A_h y = \partial \phi_h (y) = \partial \phi(y)= D\phi(y) = Ay \quad \forall y \in Y_h.
\end{equation}
As for $\psi_h$ we are allowing some extra freedom (let us however remark that \eqref{phi_psi} follows from \eqref{coerc2} as soon as $\psi_h$ is the restriction of $\psi$ to $Y_h$ since, in this case, $C_h=C \cap Y_h$). On the other hand, we are asking $\psi_h$ to be positively $1$-homogeneous, namely we are considering the case of some rate-independent approximation of \eqref{eq} only. The reader is referred instead to {\sc Efendiev \& Mielke} \cite{Efendiev-Mielke06}, {\sc Efendiev, Mielke, Rossi, \& Savar\'e} \cite{Mielke06}, and {\sc Zanini} \cite{Zanini06} for some results in the direction of rate-dependent approximation of rate-independent processes.

Finally, we shall (re)define the approximating functionals as $F_h:W^{1,1}(0,T;Y) \to [0,\infty]$ as
\begin{gather}
F_h(y)=\int_0^T \Big( \psi_h(y') + \psi^*_h(\ell_h - A_h y) - \lan \ell_h - A_h y, y'\ran \Big) + \chi_h(y(0) -y_{0,h}),\nonumber
\end{gather}
where $A_h = \partial\phi_h$ and $\chi_h(\cdot)=\phi_h(\cdot)+|\cdot|^2$. We have the following.

\begin{theorem}[Convergence of space approximations]\label{space_approx}
Assume \eqref{Y_h}-\eqref{y_0h} and let $F_h(y_h) =0$. Then $y_h \to y$ weakly star in $W^{1,\infty}(0,T;Y)$ and $F(y)=0$.
\end{theorem}
 
\begin{proof}
  By Lemma \ref{lipschitz}, we find a (not relabeled) subsequence $y_n \to y$ weakly star in $W^{1,\infty}(0,T;Y)$ and weakly pointwise. Since $F_h(y_h) = 0 $ we readily check that $y(t)\in Y_h$ for all $t \in [0,T]$. In particular, $A_h y_h = A y_h$  for all $t \in [0,T]$ owing to \eqref{identity_A}. 
Hence, by lower semicontinuity,
\begin{eqnarray}
    0 \leq F(y) &\leq& \liminf_{h \to 0}\Bigg(\int_0^T \Big( \psi_h(y'_h) + \psi^*_h(\ell_h - Ay_h) - \lan \ell_h, y'_h\ran \Big)\nonumber\\
&+& \phi(y_h(T))+ \phi(y_{0,h}) - \lan A y_h(0),y_{0,h}\ran + |y_h(0)- y_{0,h}|^2\Bigg)\nonumber\\&=&\liminf_{h \to 0}\Bigg(\int_0^T \Big( \psi_h(y'_h) + \psi^*_h(\ell_h - A_h y_h) - \lan \ell_h, y'_h\ran \Big)\nonumber\\
&+& \phi_h (y_h(T)) - \phi_h(y_{h}(0))\Bigg)\nonumber\\
&=& \liminf_{h \to 0} F_h(y_h) =0,\nonumber
  \end{eqnarray}
Note that the integral terms containing $\psi$ and $\psi^*$ pass to the $\liminf$ by means of Lemma \ref{cor44}. 
\end{proof}

 By inspecting the proof of Theorem \ref{space_approx} (which of course generalizes Lemma \ref{data_approx}), one realizes that, whenever the weak-star precompactness in $W^{1,\infty}(0,T;Y)$ of the sequence $y_h$ is assumed, the convergence statement holds more generally in the case $F_h(y_h) \to 0$. Namely, by directly asking for the above-mentioned compactness, one could consider the convergence of some approximated solutions $y_h$ such that, possibly, $F_h(y_h)>0$. We rephrase this fact in the following statement.

  \begin{lemma}[$\Gamma$-liminf inequality for $F_h$]
Assume \eqref{Y_h}-\eqref{psi_h}, \eqref{psi_h3}-\eqref{ell_h}, and \eqref{y_0h}. Then,
$$F(u) \leq \inf \left\{\liminf\limits_{h\to 0}F_h(y_h) \ : \ y_h \to y \quad \text{weakly star in} \ \ W^{1,\infty}(0,T;Y) \right\}.$$
\end{lemma}

Note that the homogeneity of $\psi_h$, the uniform convexity of $\phi_h$, and the Lipschitz continuity of $\ell_h$ play no role here. 

Finally, again by looking carefully to the proof of Theorem \ref{space_approx} one could wonder if the requirement on the Mosco convergence of $\psi_h$ could be weakened. Indeed, what we are actually using is only that 
\begin{equation}
\psi \leq \Gamma-\liminf\limits_{h\to 0}\psi_h \quad \text{and} \quad \psi^* \leq \Gamma-\liminf\limits_{h\to 0}\psi_h^*\label{gamma}
\end{equation}
with respect to the weak topologies of $Y$ and $Y^*$, respectively. On the other hand, in our specific situation, \cite[Lemma 4.1]{be} entails that \eqref{gamma} and the fact that $\psi_h \to \psi$ in the Mosco sense in $Y$ are equivalent.  

This observation motivates once again the belief that Mosco convergence is the right frame in order to pass to limits within rate-independent problems. For the sake of completeness, let us recall that a first result in the direction of the approximation of the play operator ($Y$ Hilbert and $A$ coercive on $Y$) under the Hausdorff convergence of the characteristic sets $C_h^* = D(\psi^*_h)$ is contained in \cite[Thm. 3.12, p. 34]{krejci} whereas the extension of this result to the more general situation of Mosco converging sets as well as some application to parabolic PDEs with hysteresis is discussed in \cite{lombardo}. More recently, {\sc Mielke, Roub\' \i\v cek, \& Stefanelli} \cite{mrs} addressed in full generality the issue of $\Gamma$-convergence and relaxation for the energetic solutions of rate-independent processes. An alternative convergence proof in the specific case of convex energies is obtained by means of the Brezis-Ekeland approach in \cite{be}.

\section{Time-discretization}\label{time-discretization}
\setcounter{equation}{0}

Assume now we are given the partitions $P_n=\{0=t^0_n<t^1_n<\dots<t^{N_n}_n=T\}$  and denote by $\tau^i_n= t^i_n - t^{i-1}_n$ the $i$-th time-step and by $\tau_n = \max_{1 \leq i \leq N_n}\, \tau_n^i $ the diameter of the $ n$-th partition. No constraints are imposed on the possible choice of the time-steps throughout this analysis besides $ \tau_n \to 0 $ as $ n \to \infty$. Moreover, let the parameter $\theta\in [1/2,1]$ be given. 

In the following we will make an extensive use of the following notation: letting $ v=(v^0, \dots,v^{N_n})$ be a vector, we will denote by $ \haz v_n$ and $ \ove v_n $ two functions of the time-interval $ [0,T] $ which interpolate the values of the vector $ v $ piecewise linearly and backward constantly on the partition $P_n $, respectively. Namely
\begin{gather}
  \haz v_n(0)= v^0, \ \ \haz v_n(t) = \gamma^i_n(t) v^i + \big( 1 - \gamma^i_n(t) \big) v^{i-1}, \nonumber \\
  \ove v_n(0)= v^0, \ \ \ove v_n(t)= v^i, \quad \text{for} \ \ t \in (t^{i-1}_n,t^i_n], \ \ i=1, \dots, N_n \nonumber
\end{gather}
where 
$$
\gamma^i_n(t) =(t-t^{i-1}_n)/\tau^i_n \quad \mbox{for} \ \  t \in (t^{i-1}_n,t^i_n], \  \ i=1, \dots, N_n.
$$
Moreover, we let $\delta v^i= (v^i - v^{i-1})/\tau^i_n\,$ for $\,i =1, \dots, N_n$ (so that $\haz v_n' = \overline{\delta v_n}$) and denote by $v_\theta$ the vector with components $v_\theta^i= \theta v^i + (1-\theta)v^{i-1}$.

Recall that $\ell \in W^{1,\infty}(0,T;Y^*)$ and $y_0 \in S(0)$.
We shall be concerned with the so-called $\theta$-scheme for problem \eqref{eq}:
\begin{eqnarray}
 && \partial \psi \left(\frac{y^i_n - y^{i-1}_n}{\tau^i_n}\right)+  A (\theta y^i_n + (1-\theta) y^{i-1}_n) \ni \ell(\theta t^i_n + (1-\theta) t^{i-1}_n) \nonumber\\
&& \text{for} \ \ i=1,\dots,N_n, \label{theta}\\
&& y^0_n=y_0.\label{theta2}
\end{eqnarray}
One usually refers to the latter as backward or implicit Euler scheme for the choice $\theta=1$ and as Crank-Nicholson scheme for $\theta=1/2$.

Owing to the above-introduced notation, the latter scheme can be equivalently rewritten as
\begin{gather}
  \partial \psi \left( y^i_n - y^{i-1}_n\right)+  A y^i_{n,\theta}  \ni \ell( t^i_{n,\theta}) \quad \text{for} \ \ i=1,\dots,N_n, \quad y^0_n=y_0.\label{theta3}
\end{gather}

Clearly, the $\theta$-scheme \eqref{theta3} is rate-independent. Namely,  no time-step appears in \eqref{theta3} and the choice of the partition affects the solution via the values of the load $\ell$ only. In this concern, our focus on variable time-steps partition could be simplified by considering proper rescaled loads $\ell$ instead. We shall however keep up with it, especially in order to underline the possibility of adapting the partition according to some a posteriori analysis (see Subsection \ref{adaptivity}).

Before moving on, let us comment that, for all $n$, the latter scheme as a unique solution. Indeed, given $y^{i-1}_n\in C$, it suffices to (uniquely) solve iteratively the incremental problem
\begin{equation}
  \label{incremental}
  y^i_n \in \argmin_{y \in Y} \Big( \theta \phi(y) - \lan\ell( t^i_{n,\theta}) - (1-\theta)A y^{i-1}_n, y \ran+  \psi(y - y^{i-1}_n) \Big).
\end{equation}
Note that, since $y^{i-1}_n\in C$, the  functional under minimization turns out to be uniformly convex. Hence, by \eqref{coerc1}, the minimum problem has a unique solution. In particular, exactly as in Lemma \ref{C} we have the following.

\begin{lemma}\label{C_dis} $y^i_n \in C$ for all $i=0,1,\dots,N_n$.
\end{lemma}

A crucial observation is that, as in the continuous case, the discrete trajectories show some sort of stability as well. 

\begin{lemma}[Stability of the discrete trajectories]\label{stab}
We have that
\begin{gather}
y^i_n \in \argmin_{y \in Y}\Big( \theta \phi(y) - \lan\ell( t^i_{n,\theta}) - (1-\theta)A y^{i-1}_n, y \ran + \psi(y - y^i_n) \Big)\nonumber\\
\text{for}  \ \  i=1,\dots,N_n.\label{eq_stab} 
\end{gather}
In particular, if $\theta=1$, then $y^i_n \in S(t^i_n)$.
\end{lemma}
\begin{proof}
  From the incremental formulation \eqref{incremental} and the triangle inequality for $\psi$ we get that, for all $y \in Y$,
\begin{gather}
 \theta \phi(y^i_n) - \lan\ell( t^i_{n,\theta}) - (1-\theta)A y^{i-1}_n, y^i_n \ran + \psi(y^i_n - y^{i-1}_n) \nonumber\\
\leq \theta \phi(y) - \lan\ell( t^i_{n,\theta}) - (1-\theta)A y^{i-1}_n, y \ran + \psi(y - y^{i-1}_n)\nonumber\\
\leq \theta \phi(y) - \lan\ell( t^i_{n,\theta}) - (1-\theta)A y^{i-1}_n, y \ran + \psi(y - y^{i}_n)+ \psi(y^i_n - y^{i-1}_n),\nonumber
\end{gather}
whence the assertion follows.
\end{proof}

Again as in the continuous case, we readily check that

\begin{equation}
  \label{poi_stab}
  \text{relation \eqref{eq_stab} holds iff} \quad  
\ell( t^i_{n,\theta}) - A y^i_{n,\theta} \in C^*.
\end{equation}

\subsection{The discrete variational principle}\label{discrete_principle}

 We shall now present a discrete version of the variational principle of Theorem \ref{characterization}.

We define $L^{\theta,i}_n(y,z): Y \times Y \to [0,\infty]$ as
\begin{gather}
L^{\theta,i}_n(y,z)= \psi\left(\frac{y- z}{\tau^i_n}\right) + \psi^*\left(\ell( t^i_{n,\theta}) - A(\theta  y + (1-\theta) z) \right) \nonumber\\
- \lan \ell( t^i_{n,\theta}) - A(\theta  y + (1-\theta) z), \frac{y- z}{\tau^i_n} \ran,\nonumber
\end{gather}
and the functionals $F_n^\theta: Y^{N_n +1} \to  [0,\infty]$ as
$$F_n^\theta(y^0_n, \dots, y^{N_n}_n)= \sum_{i=1}^{N_n} \tau^i_n L^{\theta,i}_n(y^i_n,y^{i-1}_n) + \chi(y^0_n-y_0).$$

\begin{lemma}[Discrete variational principle]\label{discrete_variational} 
$(y^0_n, \dots, y^{N_n}_n)$ solves \eqref{theta3} iff \linebreak$F_n^\theta(y^0_n, \dots, y^{N_n}_n)=0=\min F_n^\theta$. 
\end{lemma}

\begin{proof} Analogously to the continuous case, we have that, for all $i=1,\dots, N_n$,
 $$ \partial \psi \left(\delta y^i_n\right)+  A y^i_{n,\theta}  \ni \ell( t^i_{n,\theta}) \quad \text{iff}\quad L^{\theta,i}_n(y^i_n,y^{i-1}_n) =0,$$
and $y^0_n=y_0$ iff $\chi(y^0_n-y_0)=0$.
\end{proof}

Let us observe that the functional $F_n^\theta$ is convex and lower-semicontinuous. Moreover, by the homogeneity of $\psi$ (see \eqref{psi}), $F_n^\theta$ is actually independent of the time steps. In fact, we have
\begin{gather}
 F_n^\theta(y^0_n, \dots, y^{N_n}_n)= \sum_{i=1}^{N_n} \Bigg(\psi\left({y^i_n-y^{i-1}_n }\right) + \psi^*\left(\ell( t^i_{n,\theta}) - A  y^i_{n,\theta}) \right)\nonumber\\
 - \lan \ell( t^i_{n,\theta}) - A y^i_{n,\theta}, {y^i_n-y^{i-1}_n} \ran \Bigg) + \chi(y^0_n-y_0).\nonumber
\end{gather}  

The idea of dealing with time-discretizations via a discrete variational principle closely relates our analysis to the theory of so-called {\it variationals integrators}. The latter are numerical schemes stemming from the approximation of the action functional in Lagrangian Mechanics. By referring the reader to the monograph \cite{Hairer02} and the survey \cite{Marsden-West01}, we shall restrain here from giving a detailed presentation of the subject and limit ourselves to some (necessarily sketchy) considerations. Letting $ (t,y,p)\in [0,T] \times \Rz^m\times \Rz^m \mapsto \calL (t,y,p)$ denote the Lagrangian of a (finite dimensional, for simplicity) system, the Hamilton principle asserts that the actual trajectory $t \mapsto y(t) $ of the system minimizes the action functional
$$ y \mapsto \int_0^T\calL(t,y(t),y'(t))\, \dd t$$
among all curves with prescribed endpoints, thus solving the Lagrange equations
\begin{equation}
  \label{lagrange}
  \partial_{y_i} \calL - \ddt \partial_{p_i}\calL=0\quad \text{for $i=1, \dots , m$}.
\end{equation}
 Hence, a natural idea is that of deriving numerical schemes for Lagrangian mechanics by applying some quadrature procedure to the action functional, i.e. discretizing Hamilton's principle. The resulting discrete schemes show comparable performance with respect to other methods but generally enjoy some interesting extra (and often crucial) properties such as the conservation of suitable quantities \cite{Lew04}. Variational integrators have been intensively applied in finite-dimensional contexts and, more recently, to the situation of nonlinear wave equations \cite{Marsden98} and non-equilibrium elasticity \cite{Lew03}.

The present analysis may bear some resemblance to the above-mentioned theory. Indeed, the formulation of the $\theta$-scheme in the case $\theta=1/2$ stems exactly from the midpoint quadrature of the functional $F$  as
\begin{eqnarray} 
&&\int_{t^{i-1}_n}^{t^i_n} L(t,\haz y(t), \haz y'(t))\, \dd t \nonumber\\
&=& \tau^i_n L(t^i_{n,1/2},\haz y(t^i_{n,1/2}), \haz y'(t^i_{n,1/2}) ) \nonumber\\
&=& \psi\left({y^i-y^{i-1} }\right) + \psi^*\left(\ell( t^i_{n,1/2}) - A  \left(\frac{y^i + y^{i-1}}{2} \right) \right)\nonumber\\
&-& \lan \ell( t^i_{n,1/2}) -  A  \left(\frac{y^i + y^{i-1}}{2} \right), {y^i_n-y^{i-1}_n} \ran,\nonumber
\end{eqnarray}
where $\haz y$ is taken to be piecewise affine on the partition $P_n$.

On the other hand, our focus here is quite different. First of all, we are not dealing with the Hamilton principle (endpoints are not fixed) as we are not aimed at solving the Euler-Lagrange equations for $F$ (i.e., solve \eqref{lagrange}). Secondly, we are specifically interested at infinite-dimensional situations, namely PDEs. Finally, the only choice of $\theta$ which is directly related with a quadrature of $F$ is $\theta=1/2$ and we are not considering higher order schemes.

Before closing this discussion, let us mention that some $\Gamma$-convergence techniques have been recently exploited in the (finite-dimensional) frame of variational integrators by {\sc M\"uller \& Ortiz} \cite{Mueller-Ortiz03} (see also \cite{Maggi-Morini04}). In this same spirit, we are here providing $\Gamma$-convergence results in infinite dimensions instead.

\subsection{Stability of the $\theta$-scheme}\label{theta_stability}

It is known since {\sc Han \& Reddy} \cite{Han-Reddy95,Han-Reddy99} that the choice $\theta <1/2$ in \eqref{theta3} leads to an unconditionally unstable scheme and that, on the contrary, for $\theta\in [1/2,1]$ the $\theta$-scheme is stable in $H^1(0,T;Y)$ when $Y$ is a Hilbert space and the partitions are chosen to be uniform. 

Here we shall provide an alternative stability proof by taking into account the Banach frame. 

\begin{lemma}[Stability]\label{lemma_stability} Assume \eqref{coerc2} and let $\theta\in [1/2,1]$. Then, the solution to the $\theta$-scheme \eqref{theta3} fulfills
  \begin{equation}
    \|\haz y_{n,\theta}'\|_{L^{\infty}(0,T;Y)} \leq \frac{1}{\alpha}\| \ell'\|_{L^{\infty}(0,T;Y^*)}\quad \text{if}  \ \ \theta=1 \ \ \text{or}  \ \ \theta= \frac12 \label{dis_stab}.
\end{equation}
Moreover, for constant time-steps,
  \begin{equation} 
  \|\haz y_{n,\theta}'\|_{L^{\infty}(0,T;Y)} \leq \frac{1}{\alpha(2\theta -1)}\| \ell'\|_{L^{\infty}(0,T;Y^*)}\quad \text{if}  \ \ \frac12 < \theta <1 \label{dis_stab2}.
  \end{equation}
\end{lemma}

Our argument coincides with that of \cite[Thm. 4.4]{Mielke-Theil04} in the case of Euler, i.e. $\theta=1$ and it is an extension of the latter for the case $1/2 < \theta <1$. Here, we do not play with the variational inequality by choosing suitable tests but use the scalar relations $L^{\theta,i}_n(y^i_n,y^{i-1}_n)=0$ instead (this makes however no substantial difference since the latter scalar relations are exactly the outcome of the test on the variational inequality in  \cite[Thm. 4.4]{Mielke-Theil04}).

The stability proof for the Crank-Nicholson scheme $\theta =1/2$ is quite different from  former arguments and stems as a direct outcome of our variational approach. Let us mention that, unlike the classical parabolic situation, here the Crank-Nicholson scheme is indeed unconditionally stable. The reason for this fact is the rate-independence of the problem, namely the degenerate character of the evolution (no relaxation time). In both cases $\theta=1$ and $\theta=1/2$, the stability constant $1/\alpha$ is sharp (see Lemma \ref{lipschitz}).

We complement this analysis by providing the stability for the $\theta$-scheme for $1/2<\theta<1$ in the case  of constant time-steps (likely with a non-optimal, although explicit, stability constant).

\begin{proof} Let us prove the stability of the Crank-Nicholson scheme $\theta =1/2$ first. For this aim, it suffices to recall that 
\begin{gather}
0 = F^{1/2}_n(y^0_n, \dots, y^{N_n}_n)\nonumber\\
=\int_0^T \Big( \psi(\haz y'_n) + \psi^*(\haz \ell_n - A \haz y_n) - \lan \haz \ell_n -  A \haz y_n,  \haz y_n'\ran \Big) +  \chi(\haz y_n(0) -y_0).\nonumber
\end{gather}
Hence $\haz y_n$ minimizes the functional $F$ where $\ell$ is replaced by $\haz \ell_n$. The stability estimate follows from Lemma \ref{lipschitz}.

Let us now move to the case $1/2 < \theta \leq 1$. Relation \eqref{eq_stab} applied at level $i-1$ for some $i=2,\dots,N_n$ along with the choice $y=y^i_n$ entails that
  \begin{gather}
    \theta \phi (y^i_n - y^{i-1}_n) + \theta \phi(y^{i-1}_n) - \lan\ell( t^{i-1}_{n,\theta}) - (1-\theta)A y^{i-2}_n, y^{i-1}_n \ran \nonumber\\
\leq \theta \phi(y^i_n) - \lan\ell( t^{i-1}_{n,\theta}) - (1-\theta)A y^{i-2}_n, y^i_n \ran + \psi(y^i_n - y^{i-1}_n), \nonumber
  \end{gather}
where the extra-term $\theta \phi (y^i_n - y^{i-1}_n)$ is obtained from the fact that $\phi$ is quadratic. Hence, we have that 
\begin{eqnarray}
 &&\theta \phi (y^i_n - y^{i-1}_n) + \theta \phi(y^{i-1}_n) - \theta \phi(y^i_n) \nonumber\\
&&\leq \lan\ell( t^{i-1}_{n,\theta}), y^{i-1}_n - y^i_n
\ran + (1-\theta) \lan A(y^{i-2}_n - y^{i-1}_n),y^i_n - y^{i-1}_n  \ran + \psi(y^i_n - y^{i-1}_n) \nonumber\\
&&+ (1-\theta) \lan Ay^{i-1}_n,y^i_n - y^{i-1}_n  \ran\nonumber\\
&&=\lan\ell( t^{i-1}_{n,\theta}), y^{i-1}_n - y^i_n
\ran + (1-\theta) \lan A(y^{i-2}_n - y^{i-1}_n),y^i_n - y^{i-1}_n  \ran + \psi(y^i_n - y^{i-1}_n) \nonumber\\
&&- (1-\theta) \Big(\phi(y^{i-1}_n) + \phi(y^i_n - y^{i-1}_n ) - \phi(y^i_n) \Big),\nonumber
\end{eqnarray}
so that
\begin{gather}
  \phi(e^i_n) + \phi(y^{i-1}_n)- \phi(y^{i}_n) \leq - \lan \ell( t^{i-1}_{n,\theta}), e^i_n \ran + (\theta -1) \lan A e^{i-1}_n,e^i_n \ran + \psi(e^i_n),\label{da_sommare}
\end{gather}
where we have used $e^i_n = y^i_n - y^{i-1}_n$ in order to shorten notations.

Next, from $L^{\theta,i}_n(y^i_n,y^{i-1}_n)=0$ for $i=1,\dots,N_n$, we obtain that
\begin{eqnarray}
  0 &=& \psi(e^i_n) - \lan \ell( t^{i}_{n,\theta}) - Ay^i_{n,\theta}, e^i_n \ran\nonumber\\
&=& \psi(e^i_n) - \lan \ell( t^{i}_{n,\theta}), e^i_n\ran + \theta \Big(\phi(y^i_n) + \phi(e^i_n) - \phi(y^{i-1}_n) \Big)\nonumber\\
&-& (1- \theta) \Big( \phi(y^{i-1}_n) + \phi(e^i_n) - \phi(y^{i}_n)\Big).\nonumber
\end{eqnarray}
In particular, we have checked that
\begin{equation}\label{poppo}
\psi(e^i_n) + \phi(y^{i}_n) - \phi(y^{i-1}_n) + (2\theta -1)\phi(e^i_n) = \lan \ell( t^{i}_{n,\theta}), e^i_n\ran.
\end{equation}
We take the sum between the latter and \eqref{da_sommare} and get that
$$2\theta \phi(e^i_n) \leq \lan\ell( t^{i}_{n,\theta})- \ell( t^{i-1}_{n,\theta}), e^i_n \ran + (\theta -1) \lan A e^{i-1}_n,e^i_n\ran,$$
or, equivalently,
\begin{equation}\label{ateta}
\lan A e^i_{n,\theta}, e^i_n\ran \leq \lan\ell( t^{i}_{n,\theta})- \ell( t^{i-1}_{n,\theta}), e^i_n \ran.
\end{equation}

Now, if $\theta =1$, we conclude that
$$|e^i_n| \leq |\leq \frac{1}{\alpha}|\ell( t^{i}_{n,\theta})- \ell( t^{i-1}_{n,\theta})|_*,$$
and the assertion follows. 

In case $1/2\leq \theta <1$ and for a constant time-step partition, one proceeds from \eqref{ateta} by computing
\begin{eqnarray}
    &&\tau_n \| \haz \ell'_n\|_{L^\infty(0,T;Y^*)} \sqrt{\frac{2}{\alpha}} \sqrt{\phi(e^i_n)} \nonumber\\
&&\geq \lan\ell( t^{i}_{n,\theta})- \ell( t^{i-1}_{n,\theta}), e^i_n \ran \nonumber\\
&&\geq \theta \lan A e^i_n, e^i_n \ran + (1-\theta) \lan Ae^{i-1}_n,e^i_n \ran\nonumber\\
&&= (2\theta -1)\lan Ae^i_n,e^i_n \ran +  (1-\theta) \lan A (e^i_n +e^{i-1}_n), e^i_n\ran \nonumber\\
&&\geq 2(2\theta -1) \phi(e^i_n) + (1-\theta) \Big( \phi(e^i_n) - \phi(e^{i-1}_n)\Big)\nonumber\label{induction}
\end{eqnarray}
Note that the coefficient $(2\theta -1)$ is strictly positive as $\theta> 1/2$. By using the fact that $y^0_n = y_0 \in S(0)$ (recall \eqref{data2}), we readily check that
$$\phi(e^1_n)+ \phi(y^0_n) - \lan \ell(0), y^0_n \ran \leq \phi(y^1_n) - \lan \ell(0), y^1_n \ran + \psi(e^1_n)  $$
and, by adding the latter to \eqref{poppo} for $i=1$ we have
\begin{equation}\label{starter}
2\theta \phi(e^1_n) \leq \lan \ell(t^1_{n,\theta}) - \ell(0), e^1_n \ran \leq \tau_n\| \haz \ell'_n\|_{L^\infty(0,T;Y^*)}\sqrt{\frac{2}{\alpha}}\sqrt{\phi(e^1_n)}.
\end{equation}
Let us define
\begin{gather}
a_i^2=\phi\left(\frac{y^i_n - y^{i-1}}{\tau_n} \right)=\phi(e^i_n)/\tau_n^2, \nonumber\\
C_0= \frac{2(2\theta -1)}{1-\theta}, \ \  C_1= \frac{1}{1-\theta}\sqrt{\frac{2}{\alpha}}\| \haz \ell'_n\|_{L^\infty(0,T;Y^*)}, \ \ C_2 =\frac{C_1}{C_0}, \ \  \nonumber
\end{gather}
so that, owing to \eqref{induction}, \eqref{starter}, and using the fact that $2(2\theta -1) < 2\theta$,
\begin{eqnarray}
\big(C_0 + 1 \big) a_i^2 - a_{i-1}^2 &\leq& C_1 a_i \quad \text{for} \ \  i=2, \dots, N_n,\nonumber\\
 a_1 &\leq& \frac{1}{2\theta}\sqrt{\frac{2}{\alpha}} \| \haz \ell'_n\|_{L^\infty(0,T;Y^*)} \leq \frac{1}{2(2\theta-1)}\sqrt{\frac{2}{\alpha}} \| \haz \ell'_n\|_{L^\infty(0,T;Y^*)} \nonumber\\
&=&\frac{1}{1-\theta} \sqrt{\frac{2}{\alpha}}\| \haz \ell'_n\|_{L^\infty(0,T;Y^*)} \, \frac{1-\theta}{2(2\theta-1)}\nonumber\\
&=& \frac{C_1}{C_0} = C_2.\nonumber
\end{eqnarray}
Now, since $(C_0 +1) C_2^2 - C_1 C_2= C_2^2$, we easily prove by induction that $a_i \leq C_2$ and the assertion follows.
\end{proof}
\subsection{Convergence}\label{convergence}
 
We shall prove the weak-star $W^{1,\infty}(0,T;Y^*)$ convergence for the $\theta$-method. This result has to be compared with that of Han \& Reddy \cite[Thm. 3.4]{Han-Reddy00} where the uniform convergence of the backward constant interpolations is obtained. Our result is weaker than that of \cite[Thm. 3.4]{Han-Reddy00} since we are not providing strong convergence. On the other hand,  we believe our half-page proof to be possibly more transparent than the long argument developed in \cite{Han-Reddy00}. Let us moreover mention that in the Hilbertian case and for $A$ coercive on $Y$, the strong convergence in $W^{1,p}(0,T;Y)$ for all $p <\infty$ of the Euler method $\theta=1$ has been proved in \cite[Prop. 3.9, p. 33]{krejci}.

\begin{theorem}[Convergence for the $\theta$-method]\label{convergence_theta} 
Assume \eqref{coerc2} and let \linebreak$F_n^\theta(y^0_n, \dots, y^{N_n}_n)=0$. Then, $\haz y_n \to y$ weakly star in $W^{1,\infty}(0,T;Y)$ where \linebreak$F(y)=0=\min F$.  
\end{theorem}

\begin{proof}
  Owing to Lemma \ref{lemma_stability}, we can extract a (not relabeled) subsequence such that $\haz y_n \to y$ weakly star in $W^{1,\infty}(0,T;Y)$, hence weakly pointwise in $Y$. Moreover, we clearly have that both $\ove y_n$ and $\ove y_{n,\theta}$ converge at the same limit weakly star in $L^\infty(0,T;Y)$. Finally, we directly check that $\ove \ell_{n,\theta} \to \ell$ strongly in $L^\infty(0,T;Y^*)$. By observing that, since $\theta \geq 1/2$,
  \begin{eqnarray}
   \tau^i_n \lan A(\theta y^i_n + (1-\theta) y^{i-1}_n), \delta y^i_n \ran &=& \phi(y^i_n)+ (2\theta -1) \phi(y^i_n - y^{i-1}_n) - \phi(y^{i-1}_n) \nonumber\\
 &\geq& \phi(y^i_n) - \phi(y^{i-1}_n),\nonumber
  \end{eqnarray}
we compute that
  \begin{eqnarray}
    0 & =&F_n^\theta(y^0_n, \dots, y^{N_n}_n) \nonumber\\
&\geq&  \int_0^T \Big(\psi(\haz y'_n) + \psi^* (\ove \ell_{n,\theta} - A \ove y_{n,\theta}) - \lan \ove \ell_{n,\theta}, \haz y'_n\ran \Big) \nonumber\\
&+& \phi(\haz y_n(T)) - \phi(\haz y_n(0)) + \chi(\haz y_n(0) - y_0)\nonumber\\
&=& \int_0^T \Big(\psi(\haz y'_n) + \psi^* (\ove \ell_{n,\theta} - A \ove y_{n,\theta}) - \lan \ove \ell_{n,\theta}, \haz y'_n\ran \Big) \nonumber\\
&+& \phi(\haz y_n(T)) + \phi(y_0) - \lan A\haz y_n(0) , y_0\ran + |\haz y_n(0) - y_0|^2.\nonumber
  \end{eqnarray}
Finally, it suffices to pass to the $\liminf$ above as $n \to \infty$ and exploit lower semicontinuity and the stated convergences in order to obtain that $F(y)\leq 0$. Hence, by Theorem \ref{characterization} and Corollary \ref{uniqueness}, $y$ is the only solution to \eqref{eq} and the whole sequence $\haz y_n$ converges. 
\end{proof}

We shall mention that the separability assumption for $Y$ is not crucial and could be weakened. Indeed, in case $Y$ is not separable, one simply has to pass to limits as $\haz y_n \to y$ weakly in $W^{1,p}(0,T;Y)$ for some $p\in [1,\infty)$.

\subsection{The functional controls the uniform distance}\label{distance_euler}

We shall reproduce at the discrete level the results of Subsection \ref{uniqueness}. We begin by showing how to possibly control the uniform distance of two vectors by means of the discrete functional $F_n^\theta$.  

\begin{lemma}[Uniform distance control via $F_n^\theta$]\label{uniform_distance_dis} Let the vectors $u=(u^0, \dots, u^{N_n})$,  $ v=(v^0, \dots, v^{N_n}) \in Y^{N_n+1}$ be given. Then,
\begin{gather}
\eta(1-\eta) \max_{1\leq i \leq N_n}\phi(u^i - v^i) \nonumber\\
\leq\eta F_n^\theta(u^0, \dots, u^{N_n}) + (1-\eta)F_n^\theta(v^0, \dots, v^{N_n})\quad
 \forall \eta \in [0,1].\nonumber
\end{gather}
\end{lemma}

\begin{proof}
This proof follows the same lines of that of Corollary \ref{uniform_distance}. Let $ 1 \leq i  \leq N_n$ be fixed and define $F^{\theta, i}_n: Y^{i+1} \to [0,\infty]$ as
\begin{eqnarray}
&&\!\!\!\!\!\!\!\!\!\!\!\!\!\!\!\!F^{\theta,i}_n(y^0, \dots, y^i)\nonumber\\
&=& \sum_{j=1}^i \tau^j_n L^{\theta,j}_n(y^j,y^{j-1}) + \chi(y^0 - y_0)\nonumber\\
&=& \sum_{j=1}^i \psi(y^j - y^{j-1}) + \psi^*(\ell(t^j_{n,\theta}) - A y^j_\theta) - \lan\ell(t^j_{n,\theta}),y^j - y^{j-1} \ran  \nonumber\\
&+& \phi(y^i)  + (2 \theta -1 )\sum_{j=1}^i \phi(y^j - y^{j-1})+ \phi(y_0)  - \lan A y^0, y_0\ran + |y^0-y_0|^2. \nonumber
\end{eqnarray}
Then, clearly $y=(y^0, \dots, y^i)\mapsto G^{\theta,i}_n(y)=F^{\theta,i}_n(y) - \phi(y^i)$ is convex. Hence, letting $w=\eta u + (1-\eta)v$ for some $\eta \in [0,1]$ we have that 
\begin{gather}
  0 \leq F^{\theta,i}_n(w) \nonumber\\
\leq \eta \big(G^{\theta,i}_n(u) + \phi(u^i)\big) + (1-\eta)\big(G^{\theta,i}_n(v) + \phi(v^i)\big) - \eta(1-\eta)\phi(u^i- v^i), \nonumber
\end{gather}
whence the assertion follows.
\end{proof}

Again, note that the latter lemma controls the uniform norm of the distance only if the stronger \eqref{coerc2} is required. The following corollary of Lemma \ref{uniform_distance_dis} will be the starting point for some possible a posteriori error control procedure (see Subsection \ref{a_posteriori}).

\begin{corollary}[Uniform distance from the minimizer]\label{unif_dist_min_dis}
Let $F_n^\theta(y^0, \dots, y^{N_n})=0$. Then
\begin{gather}
\max_{1\leq i \leq N_n}\phi(y^i- v^i) \leq F_n^\theta(v^0, \dots, v^{N_n}) \quad \forall (v^0, \dots, v^{N_n}) \in Y^{N_n +1}.\nonumber
\end{gather}
\end{corollary}

Moreover, we re-obtain a proof of the uniqueness of the solution of the $\theta$-method.

\begin{corollary}[Uniqueness of the minimizer]\label{un_dis} 
Assume \eqref{coerc2}. Then, there exists at most one  $y=(y^0, \dots, y^{N_n})$ such that $F_n^\theta(y)=0$.
\end{corollary}

\subsection{The generalized $\theta$-method}\label{discrete_scheme}

Although minimizers of $F^\theta_n$ and solutions of the $\theta$-scheme \eqref{theta3} coincide, minimizing sequences of $F^\theta_n$ need not solve \eqref{theta3}. This extra freedom allows the minimization formulation to capture the convergence of some generalized $\theta$-method, where the relations in \eqref{theta3} are not solved exactly but rather approximated. Namely, we shall look for vectors $u_n=(u^0_n, \dots, u^{N_n}_n)$ such that
$$F^\theta_n (u_n) \to 0 \quad \text{as} \quad n \to \infty$$
instead of $F^\theta_n (u_n)=0$ for all $n\in \Nz$.

From the computational viewpoint, note that the $\theta$-scheme consists in solving $N_n$ nonlinear equations in one unknown each while checking for stationarity for $F^\theta_n$ implies the solution of a tridiagonal system of $\, N_n+1 \,$ nonlinear equations with (up to) three unknowns each. This entails in particular that minimizing $F^\theta_n$ instead of solving \eqref{theta3} could be of a scarce interest if one is merely concerned in reproducing the $\theta$-scheme with no error. On the other hand, the issue of solving up to some tolerance turns out to be particularly relevant whenever one is aimed at implementing an optimization procedure for the solution of \eqref{theta3}. Indeed, one should be prepared to run the algorithm (some descent method, say) until some given tolerance is reached. 

Our starting point for a possible convergence analysis of the generalized $\theta$-method is the following classical error control result.

\begin{theorem}[Mielke \& Theil \cite{Mielke-Theil04}] Assume \eqref{coerc2}. Then, $\haz y_n \to y$ uniformly and $F(y)=0$. In particular,
\begin{equation}
\max_{t \in [0,T]}|( \haz y_n - y )(t)| \leq C_{e} \sqrt{\tau_n}\label{error_euler} 
\end{equation}
where $C_{e}$ depends only on data and is independent of $n$.
\end{theorem}

More precisely, in \cite{Mielke-Theil04} solely the case of the Euler scheme $\theta=1$ is discussed. However, an easy adaptation of the argument entails the result for $\theta \in [1/2,1)$ as well.

By explicitly comparing the minimizing sequence $u_n=(u^0_n, \dots, u^{N_n}_n)$ with the corresponding solution $(y^0_n, \dots,y^{N_n}_n)$ of the $\theta$-method, we have the following.

\begin{theorem}[Convergence for the generalized $\theta$-method]\label{conv_theta}
  Assume \eqref{coerc2} and let \linebreak$F_n^\theta(u^0_n, \dots, u^{N_n}_n) \to 0$. Then, $\haz u_n \to y$ uniformly, where $F(y)=0$. In particular,
\begin{equation}
\max_{t \in [0,T]}|( \haz u_n - y )(t)| \leq C_e \sqrt{\tau_n} + \left(\frac{2}{\alpha} F_{n }^\theta(u^0_n, \dots, u^{N_n}_n)\right)^{1/2}.\label{error_euler_2} 
\end{equation}
\end{theorem}
\begin{proof}
  We have that
\begin{eqnarray}
\max_{t \in [0,T]}|(y - \haz u_n  )(t)| &\leq& \max_{t \in [0,T]}|(y -  \haz y_n  )(t)| + \max_{t \in [0,T]}|( \haz y_n - \haz u_n )(t)| \nonumber\\
&\leq&  C_e \sqrt{\tau_n} + \max_{1\leq i \leq N_n}|u^i_n - y^i_n| \nonumber\\
&\leq& C_e \sqrt{\tau_n} + \left(\frac{2}{\alpha} \max_{1\leq i \leq N_n}\phi(u^i_n - y^i_n)\right)^{1/2}\nonumber
\end{eqnarray}
and we conclude by applying Corollary \ref{unif_dist_min_dis}.
\end{proof}


\subsection{A posteriori error control}\label{a_posteriori}
Let us now exploit both Corollary \ref{unif_dist} and Theorem \ref{conv_theta} in order to provide some possible a posteriori estimates of the approximation error by means of solutions $u_n$ of the generalized $\theta$-method described above.

\begin{lemma}[A posteriori error control via $F_n^\theta$] Assume \eqref{coerc2} and let \linebreak $F_n^\theta(u^0_n, \dots, u^{N_n}_n) \sim \tau_n^s$ for some $s >0$ and $F(y)=0$. Then, 
$$\max_{t \in [0,T]}|( \haz u_n - y )(t)| \sim \tau_n^r \quad \text{where}\quad 2r=\max\{1,s\}.$$
\end{lemma}

\begin{lemma}[A posteriori error control via $F$] Assume \eqref{coerc2} and let $F(\haz u_n)\sim \tau_n^s$ for some $s >0$ and $F(y)=0$. Then, $\max_{t \in [0,T]}|( \haz u_n - y )(t)| \sim \tau_n^{s/2}$.
\end{lemma}

We are also in the position of proving the weak-star convergence of the time derivatives of solutions $u_n$ of the generalized $\theta$-method by comparing them with the corresponding derivatives of the exact solution of the $\theta$-method. 
\begin{lemma}[Improved convergence for the generalized $\theta$-method] Assume \eqref{coerc2} and let $F^\theta_n(u_n^0, \dots, u^{N_n}_n )\sim \tau_n^2$. Then, $\haz u_n$ is equibounded in $W^{1,\infty}(0,T;Y)$. In particular, $\haz u_n \to y$ weakly star in $W^{1,\infty}(0,T;Y)$.
\end{lemma}

\begin{proof}
  Let $(y^0_n, \dots, y^{N_n}_n)$ be the solution of the $\theta$-scheme. By exploiting Lemma \ref{unif_dist_min_dis}, we check that
  \begin{eqnarray}
     | u^i_n - u^{i-1}_n|
&\leq& \left|{u^i_n - y^i_n} \right|+ \left|{y^i_n - y^{i-1}_n}\right| + \left|y^{i-1}_n - u^{i-1}_n \right| \nonumber\\
&\leq & \tau^i_n\|\haz y_n'\|_{L^{\infty}(0,T;Y)} + 2\left(\frac{2}{\alpha} F_n^\theta(u^0_n, \dots, u^{N_n}_n) \right)^{1/2}.\nonumber
  \end{eqnarray}
The uniform bound on $\| \haz u_n \|_{W^{1,\infty}(0,T;Y)}$ follows by dividing the latter by $\tau^i_n$, taking the maximum as $1\leq i \leq N_n$, and recalling Lemma \ref{lemma_stability}. 
\end{proof}

\subsection{Adaptivity}\label{adaptivity}

Assuming \eqref{coerc2}, the above introduced a posteriori error estimators can be exploited in order to develop an adaptive strategy. In particular, the error control in the uniform norm up to a given tolerance $\text{\rm tol}>0$ 
$$\max_{t\in [0,T]}|(y - \haz y_n)(t)| \leq  \text{\rm tol}$$
for some piecewise approximation $\haz y_n$ with $\chi(\haz y_n(0) - y_0)\leq \alpha \, \text{\rm tol}^2/4$ can be inferred, for instance, by choosing time steps in such a way that 
$$ \int_{t^{i-1}_n}^{t^i_n} L(t,\haz y_n(t),\haz y_n'(t)) \leq \frac{\alpha\, \text{\rm tol}^2}{4 N_n}.$$
Namely, by uniformly distributing the error along the partition.

Alternatively, one could develop an adaptive strategy by considering just computed quantities at the discrete level by asking for 
$$ \tau^i_n L^{\theta,i}_n(y^i_n,y^{i-1}_n) \leq \frac{\alpha\, \text{\rm tol}^2}{32 N_n} \quad \text{for} \quad \tau_n \leq \frac{\text{\rm tol}^2}{ 16 C_e^2}$$
and exploiting Theorem \ref{conv_theta}.

\section{Space-time approximations}\label{space-time-approximations}
\setcounter{equation}{0}

Let us combine the results of the previous sections (and use the corresponding notation) in order to state and prove a result on the convergence of full space-time approximations. Our results have to be compared with the former convergence analysis by {\sc Han \& Reddy} \cite{Han-Reddy99}. Our approach leads to a convergence proof with respect to a weaker topology. However, it is on the one hand slightly more general (some assumptions on the spaces and the functionals, see (H1)-(H2) \cite[p. 264]{Han-Reddy99}, are not required) and on the other hand has a much simpler proof.

\begin{theorem}[Convergence of space-time approximations]\label{s_t}
Assume \eqref{Y_h}-\eqref{y_0h},\linebreak $\theta \in [1/2,1]$, define $L^{\theta,i}_{n,h}(y,z): Y \times Y \to [0,\infty]$ as
\begin{gather}
L^{\theta,i}_{n,h}(y,z)= \psi_h\left(\frac{y- z}{\tau^i_n}\right) + \psi^*_h\left(\ell_h( t^i_{n,\theta}) - A_h(\theta  y + (1-\theta) z) \right) \nonumber\\
- \lan \ell_h( t^i_{n,\theta}) - A_h(\theta  y + (1-\theta) z), \frac{y- z}{\tau^i_n} \ran,\nonumber
\end{gather}
where $A_h = \partial\phi_h$, and let the functionals $F^\theta_{n,h}: Y^{N_n +1} \to  [0,\infty]$ be defined as
$$F^\theta_{n,h}(y^0, \dots, y^{N_n})= \sum_{i=1}^{N_n} \tau^i_n L^{\theta,i}_{n,h}(y^i,y^{i-1}) + \chi_h(y^0-y_0),$$
where $\chi_h(\cdot)=\phi_h(\cdot)+|\cdot|^2$
(note that $D(F_{n,h}^\theta) \subset Y^{N_n +1}_h $).
Finally, let \linebreak$F_{n,h}(y^0_h, \dots, y^{N_n}_h)=0$. We have:
\begin{itemize}
\item[(a)] $\haz y_{n,h} \to y_h$ weakly star in $W^{1,\infty}(0,T;Y)$ as $(n,h)\to (\infty,h)$ and $F_h(y_h)=0$.
\item[(b)] $\haz y_{n,h} \to \haz y_n$ weakly star in $W^{1,\infty}(0,T;Y)$ as $(n,h)\to (n,0)$ and $F^\theta_n(y_n)=0$.
\item[(c)] $\haz y_{n} \to y$ weakly star in $W^{1,\infty}(0,T;Y)$ as $(n,0)\to (\infty,0)$ and $F(y)=0$.
\item[(d)]$ y_{h} \to y$ weakly star $W^{1,\infty}(0,T;Y)$ as $(\infty,h)\to (\infty,0)$ and $F(y)=0$.
\item[(e)] $\haz y_{n,h} \to y$ weakly star in $W^{1,\infty}(0,T;Y)$ as $(n,h)\to (\infty,0)$ and $F(y)=0$.
\end{itemize}
\end{theorem}

The thesis of the Theorem is illustrated in Figure 1. In particular, we aim at showing that the space (or data) and time-limit can be taken in any order. Note that Limit (c) has been already checked in Theorem \ref{convergence_theta} and that the very same argument yields Limit (a) as well (recall that $Y_h$ is closed). Moreover, Limit (d) is discussed in Theorem \ref{space_approx}. So what we are actually left to check are Limits (b) and (e) only.

\begin{figure}[htbp]
  \centering
\psfrag{r}{$n$}
\psfrag{t}{$h$}
\psfrag{0}{$(\infty,0)$}
\psfrag{a}{$a$}
\psfrag{b}{$b$}   
\psfrag{c}{$c$}
\psfrag{d}{$d$}
\psfrag{e}{$e$} 
 \includegraphics[width=.4 \textwidth]{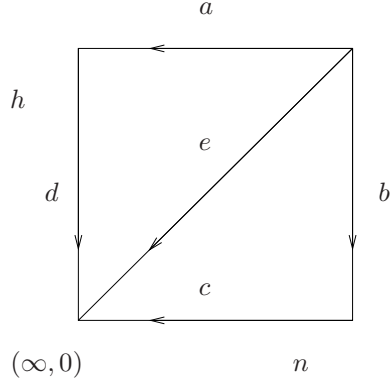}
  \caption{Convergences for space-time approximations}
  \label{fig:1}
\end{figure}

\begin{proof}
Limit (b). The assertion follows once we check that, for all $i=1,  \dots, N_n$, if $y^{i-1}_{n,h} \to y^{i-1}_{n}$ weakly in $Y$, one has the weak convergence $y^{i}_{n,h} \to y^{i}_{n}$ as well. Recall that 
\begin{gather}
  y^i_{n,h} \in \argmin_{y \in Y} \Big( \theta \phi_h(y) - \lan\ell_h( t^i_{n,\theta}) - (1-\theta)A_h y^{i-1}_{n,h}, y \ran+  \psi_h(y - y^{i-1}_{n,h}) \Big)\nonumber\\
= \argmin_{y \in Y_h} \Big( \theta \phi(y) - \lan\ell_h( t^i_{n,\theta}) - (1-\theta)A y^{i-1}_{n,h}, y \ran+  \psi_h(y - y^{i-1}_{n,h}) \Big).\nonumber
\end{gather}
Hence, since we have \eqref{phi_psi}, the sequence $y^i_{n,h}$ is weakly precompact and, up to the extraction of a (not relabeled) subsequence, $y^i_{n,h} \to \tilde y$ weakly in $Y$. Let us prove that $\tilde y$ solves the incremental problem \eqref{incremental}. Indeed, we have that
\begin{eqnarray}
0 &\leq& L^{\theta,i}_n(\tilde y,y^{i-1}_n) \nonumber\\
  &\leq & \liminf\limits_{h \to 0} \Big(\psi_h(y^{i}_{n,h} - y^{i-1}_{n,h}) + \psi^*_h (\ell_h( t^i_{n,\theta}) - Ay^i_{n,h,\theta}) \nonumber\\
&-& \lan\ell_h( t^i_{n,\theta}) - Ay^i_{n,h,\theta} , y^{i}_{n,h} - y^{i-1}_{n,h}\ran   \Big) \nonumber\\
&=&\liminf\limits_{h \to 0} L^{\theta,i}_{n,h}(y^{i}_{n,h}, y^{i-1}_{n,h})=0\nonumber
\end{eqnarray}
where we have used the Mosco convergence in \eqref{psi_h3} and the pointwise convergence of $\ell_h$ \eqref{ell_h}. Since the only solution of \eqref{incremental} is $y^i_n$, we have that $\tilde y= y^i_n$ and the whole sequence converges.

Let us mention that, if the functionals $\psi_h$ are uniformly linearly bounded (which is quite common in practice), then one could prove the latter convergence to be actually strong: namely $y^{i-1}_{n,h} \to y^{i-1}_n$ strongly in $Y$ implies the strong convergence $y^{i}_{n,h} \to y^{i}_n$. Indeed, let $w_h$ and $\tilde w_h$ be such that $w_h - y^i_{n,h} \to 0$ strongly in $Y$, $\psi_h(w_h - y^i_{n,h})\to 0$, $\tilde w_h \in Y_h$ and $\tilde w_h -w_h \to 0$ strongly in $Y$. Then
\begin{gather}
\theta \phi( y^i_{n,h}) - \lan\ell( t^i_{n,\theta}) - (1-\theta)A_h y^{i-1}_{n,h},  y^i_{n,h} \ran  + \psi_h(y^i_{n,h} - y^{i-1}_{n,h})\nonumber\\
\leq\theta \phi_h(\tilde w_h) - \lan\ell( t^i_{n,\theta}) - (1-\theta)A y^{i-1}_{n,h}, \tilde w_h \ran+  \psi_h(\tilde w_h - w_h) + \psi_h(w_h - y^{i}_{n,h}).\nonumber
\end{gather}
If $\psi_h$ are uniformly linearly bounded above, then $\psi_h(\tilde w_h - w_h)\to 0 $ with $h \to 0$. 
Whence, by passing to the $\limsup$ in the latter, we check that $\limsup_{h \to 0}\phi(y^i_{n,h}) \leq \phi(y^i_n)$, which, together with lower semicontinuity gives $\phi(y^i_{n,h}) \to \phi(y^i_n)$ and the strong convergence follows from the reflexivity of $Y$.

Limit (e). Lemma \ref{lemma_stability}, the uniform Lipschitz continuity of $\ell_h$ \eqref{ell_h2}, and the initial datum convergence \eqref{y_0h}, entail that $\haz y_{n,h}$ are uniformly Lipschitz continuous as well. Hence, by extracting a (not relabeled) subsequence, $\haz y_{n,h} \to y$ weakly star in $W^{1,\infty}(0,T;Y)$. In order to check that $y$ solves \eqref{eq}, let us remark that, being $\ell_{n,h,\theta}^i= \ell_h(t^i_{n,\theta})$,
$$\ove \ell_{n,h,\theta} \to \ell \quad\text{strongly in} \ \ L^1(0,T;Y^*)$$
and that, by \cite[Cor. 4.4]{be}
\begin{gather}
\int_0^T \psi(y') \leq \liminf\limits_{h\to 0}\int_0^T \psi_h(\haz y_{n,h}')\nonumber\\
\int_0^T \psi^*(\ell - A y) \leq \liminf\limits_{h\to 0} \int_0^T \psi_h^*(\ove \ell_{n,h,\theta} - A \ove y_{n,h,\theta}),\nonumber
\end{gather}
and compute that
\begin{eqnarray}
0 &\leq& F(y)\nonumber\\
&\leq& \liminf\limits_{h\to 0}\Bigg(\int_0^T \big(\psi_h(\haz y_{n,h}') + \psi_h^*(\ove \ell_{n,h,\theta} - A \ove y_{n,h,\theta})- \lan \ove \ell_{n,h,\theta} , \haz y_{n,h}'\ran \big) \nonumber\\
&+& \phi(\haz y_{n,h}(T)) + \phi(y_{0,h}) - \lan A\haz y_{n,h}(0), y_{0,h}\ran+ |\haz y_{n,h}(0) - y_{0,h}|^2  \Bigg)\nonumber\\
&\leq& \liminf\limits_{h\to 0} F_{n,h}^\theta(y^0_{n,h},\dots,y^{N_n}_{n,h}) =0,\nonumber
\end{eqnarray}
and we have that $F(y)=0$.
\end{proof}

We shall conclude by briefly mention some further results which can be obtained by suitably adapting to the current fully-discretized situation the arguments developed above for time-discretizations.
Firstly, in the same spirit of Lemma \ref{uniform_distance_dis},
one could consider the possibility of estimating the distance  of a vector from the minimizer of  $F^\theta_{n,h}$ by means of the functional itself.
Secondly, the use of Corollary \ref{unif_dist} would entail the possibility of an a posteriori error control an some adaptive strategy along the lines of Subsection \ref{adaptivity} could be considered.
Finally, by relying on the known convergence estimates for full space-time discretized problems \cite{Han-Reddy99} one could obtain a convergence and an a posteriori error control result for some generalized space-time approximated problem where $F^\theta_{n,h}$ are not exactly minimized and one considers minimizing sequences instead (see Subsection \ref{discrete_scheme}). We shall develop these considerations elsewhere.

\bibliographystyle{alpha}
\def\cprime{$'$} \def\cprime{$'$}

\end{document}